\documentclass[12pt]{amsart}

\usepackage{hyperref}
\usepackage{amssymb,amsmath}
\usepackage[alphabetic,nobysame,abbrev]{amsrefs}
\usepackage{enumerate}
\usepackage{color}
\usepackage{pdfsync}
\RequirePackage{fix-cm}

\DeclareFontFamily{U}{mathx}{}
\DeclareFontShape{U}{mathx}{m}{n}{<-> mathx10}{}
\DeclareSymbolFont{mathx}{U}{mathx}{m}{n}
\newlength\savedwidth 
\newcommand\whline{\noalign{\global\savedwidth\arrayrulewidth\global\arrayrulewidth 1.2pt}%
	\hline
	\noalign{\global\arrayrulewidth\savedwidth}}

\usepackage{fancybox,framed}
\usepackage{bm}
\usepackage{mdframed}
\usepackage{cancel}
\usepackage{makecell}
\usepackage{pdfpages}
\usepackage{pgf,tikz}
\usepackage{mathrsfs}
\usetikzlibrary{arrows}
\usetikzlibrary{arrows,decorations.markings}
\usetikzlibrary{matrix}
\usetikzlibrary{backgrounds}

	\usepackage{diagbox}

\definecolor{qqqqff}{rgb}{0.,0.,1.}
\definecolor{xdxdff}{rgb}{0.49019607843137253,0.49019607843137253,1.}

\definecolor{qqqqff}{rgb}{0.,0.,1.}

\usepackage{tikz}
\usepackage{graphicx}
\usepackage[cmtip,all]{xy}

\usepackage[draft]{todonotes}
\usepackage[cmtip,all]{xy}

\usepackage{hyperref}
\usepackage[utf8]{inputenc}
\usepackage{amsmath,amsfonts,amssymb,euscript,graphicx,color,tikz}

\theoremstyle{plain}
\newtheorem{theorem}{Theorem}[subsection]

\newtheorem{lem}[theorem]{Lemma}

\newtheorem{pro}[theorem]{Proposition}

\theoremstyle{definition}
\newtheorem{DEF}[theorem]{Definition}

\newtheorem{rem}[theorem]{Remark}

\newtheorem{parag}[theorem]{{}}

\numberwithin{equation}{section}

\newcommand{\sub}{\subseteq}

\newcommand{\la}{Lie algebra }

\newcommand{\fm}{(\cdot,\cdot)}

\newcommand{\Z}{\mathbb{Z} }

\newcommand{\fg}{\mathfrak{g}}

\newcommand{\ep}{\hfill$\Box$}

\def\ad{\hbox{ad}}
\def\andd{\quad\hbox{and}\quad}
\def\sg{\sigma}
\def\a{\alpha}
\def\b{\beta}
\def\lam{\lambda}
\def\Lam{\Lambda}
\def\ep{\epsilon}
\def\andd{\quad\hbox{and}\quad}
\def\supp{\hbox{supp}}
\def\id{\hbox{id}}

\def\andd{\quad\hbox{and}\quad}
\def\ind{\hbox{ind}}

\def\fm{(\cdot,\cdot)}
\def\a{\alpha}

\def\sub{\subseteq}
\def\rd{\dot{R}}

\def\lam{\lambda}
\def\Lam{\Lambda}

\def\1k{\frac{1}{k}}

\def\la{\langle}
\def\ra{\rangle}
\def\rds{\dot{R}_{sh}}
\def\rdl{\dot{R}_{lg}}

\def\b{\beta}

\def\qed{\hfill$\Box$}

\def\sg{\sigma}

\def\hh{{\mathcal H}}

\def\sg{\sigma}

\usepackage{verbatim} 

\def\quadd{\quad\quad}

\def\ad{\hbox{ad}}

\def\be{{\bf e}}
\def\bq{{\bf q}}

\def\bbbz{{\mathbb Z}}

\def\bbbr{{\mathbb R}}

\def\bbbk{{\mathbb K}}

\def\cc{{\mathcal C}}
\def\dd{\mathcal D}
\def\ep{\epsilon}

\def\ll{{\mathcal G }}

\def\jj{{\mathcal J}}
\def\ll{\mathcal L}

\def\supp{\hbox{supp}}
\def\rad{\hbox{rad}
	\def\b{\beta}}
\def\bb{{\mathcal B}}

\def\proof{{\noindent\bf Proof. }}

\def\rds{\dot{R}_{sh}}
\def\rdl{\dot{R}_{lg}}

\def\cent{\hbox{Cent}}

\def\scd{\hbox{SCDer}}
\def\Hom{\hbox{Hom}}
\def\da{\dot\alpha}

\def\rank{\hbox{rank}}

\def\ss{\mathfrak{S}}
\usepackage[utf8]{inputenc}  
\usepackage{lmodern}          
\usepackage[T1]{fontenc}      

\def\DynkinNodeSize{1.5mm}
\def\DynkinArrowLength{2mm}
\tikzset{
	dnode/.style={
		circle,
		inner sep=0pt,
		minimum size=\DynkinNodeSize,
		fill=white,
		draw},
	middlearrow/.style={
		decoration={markings,
			mark=at position 0.8 with
			{\rdaw (0:0mm) -- +(+140:\DynkinArrowLength); \rdaw (0:0mm) -- +(-140:\DynkinArrowLength);},
		},
		postaction={decorate}
	},
	leftrightarrow/.style={
		decoration={markings,
			mark=at position 0.999 with
			{
				\rdaw (0:0mm) -- +(+135:\DynkinArrowLength); \rdaw (0:0mm) -- +(-135:\DynkinArrowLength);
			},
			mark=at position 0.001 with
			{
				\rdaw (0:0mm) -- +(+45:\DynkinArrowLength); \rdaw (0:0mm) -- +(-45:\DynkinArrowLength);
			},
		},
		postaction={decorate}
	},
	sedge/.style={
	},
	dedge/.style={
		middlearrow,
		double distance=0.6mm,
	},
	tedge/.style={
		middlearrow,
		double distance=1.0mm+\pgflinewidth,
		postaction={draw}, 
	},
	infedge/.style={
		leftrightarrow,
		double distance=0.5mm,
	},
}

\begin{document}

%
%
\title{Chevalley bases for elliptic extended affine Lie algebras of type $A_1$}

\author{S. Azam}
\address
{Department of Pure Mathematics\\Faculty of Mathematics and Statistics\\
	University of Isfahan\\ P.O.Box: 81746-73441\\ Isfahan, Iran, and\\
	School of Mathematics, Institute for
	Research in Fundamental Sciences (IPM), P.O. Box: 19395-5746.} \email{azam@ipm.ir, azam@sci.ui.ac.ir}
\thanks{This work is based upon research funded by Iran National Science Foundation (INSF) under project No.  4001480.}
\thanks{This research was in part carried out in
	IPM-Isfahan Branch.}


\begin{abstract}
	We investigate Chevalley bases for extended affine Lie algebras of type $A_1$. 
	The concept of integral structures for extended affine Lie algebras of rank greater than one has been successfully explored in recent years. However, for the rank one it has turned out that the situation becomes more delicate. In this work, we consider $A_1$-type extended affine Lie algebras of {nullity} $2$, known as elliptic extended affine Lie algebras. These Lie algebras are build using the Tits-Kantor-Koecher (TKK) construction by applying some specific Jordan algebras: the plus algebra of a quantum torus, the Hermitian Jordan algebra of the ring of Laurent polynomials equipped with an involution, and the Jordan algebra associated with a semilattice. By examining these ingredient we determine appropriate bases for null spaces of the corresponding elliptic extended affine Lie algebra leading to the establishment of Chevalley bases for these Lie algebras.
\end{abstract}
\maketitle

\setcounter{equation}{-1}
\section{\bf Introduction}\setcounter{equation}{0}
\label{intro}
In this work, we investigate the construction of Chevalley bases for elliptic extended affine Lie algebras (EALAs) of type $A_1$. A Chevalley basis, a concept introduced by Claude Chevalley, plays a central role in the structure theory of Lie algebras; see \cite{Che55}, \cite{Ste16}, \cite{Bou08} and \cite{Hum72}. It provides a systematic framework for constructing Lie algebras from their root systems and offers a unified approach to studying finite-dimensional simple Lie algebras. The importance of Chevalley bases extends beyond classical theory, providing the first unified construction of both classical and exceptional groups over arbitrary fields. This idea was later generalized to the theory of infinite-dimensional Kac-Moody algebras; see, for example,  \cite{Gar78},\cite{Gar80} and \cite{Mit85}. In recent years and in the context of modular theory, the concepts of Chevalley bases and integral structures have been explored for the class of EALAs; see \cite{AFI22}, \cite{AI23}.

Extended affine Lie algebras (EALAs) over the filed $\bbbk$ of complex numbers were introduced in 1990 by \cite{H-KT90} and were systematically studied in \cite{AABGP97}. These Lie algebras expand on the properties of finite-dimensional simple and affine Lie algebras. EALAs capture the symmetry structures that extend beyond the scope of classical and affine Lie algebras. The construction of an EALA involves an extended Cartan subalgebra that leads to a root space decomposition, along with central and derivational elements. A method developed by Erhard Neher \cite{Neh04},\cite{Neh11} provides a framework for constructing any (tame) extended affine Lie algebra starting from a centerless Lie torus and adding appropriate central and derivation elements. The nullity $\nu$ of an EALA is defined as the dimension of the real span of isotropic roots.

As was mentioned in the previous paragraph, Lie tori are essential in the construction of EALAs. They in fact characterize the cores of extended affine Lie algebras. More precisely, the core modulo center of an extended affine Lie algebra is a Lie torus, and conversely any centerless Lie torus is isomorphic to the core modulo center of an EALA. 
It is known that  centerless Lie tori associated to $A_1$-type EALAs can be obtained by Tits-Kantor-Koecher-construction (TKK-construction) that
assigns a Lie algebra to a Jordan algebra. 
According to \cite{Yos00}, the Jordan algebras involved in these constructions are:
\begin{itemize}
	\item \( \mathbb{K}^+_{\mathbf{q}} \), the plus algebra of quantum torus $\bbbk_\bq$ which is the associative algebra over $\bbbk$ defined by generators $x_i^{\pm1}$, $1\leq i\leq\nu$ and relations $x_i x_i^{-1}=1=x_i^{-1}x_i$, $x_i x_j=q_{ij}x_j x_i$. Here $\bq$ is a quantum matrix \( \mathbf{q} = (q_{ij}) \) with entries satisfying \( q_{ii} = 1 \) and \( q_{ij}q_{ji} = 1 \),  
	\item \( H(\mathbb{K}_{\be}, -) \), the \emph{Hermitian Jordan algebra} over \( \mathbb{K}_{\be} \) with the involution induced by \( \overline{x_i} = x_i \). Here $\be$ is the quantum matrix with all entries equal to $1$,  
	\item \( \mathcal{J}_S \), the \emph{Jordan algebra} associated with a semilattice \( S \) within a free abelian group of rank $\nu$, and  
	\item \( \mathbb{A}_t \), which corresponds to algebras of nullity $>2$ that will not be treated in this work.  
\end{itemize}

The paper is arranged as follows. Section \ref{preliminaries} presents some basic definitions and facts about extended affine Lie algebras and root systems, which we need in this work.  Section \ref{tame6} reviews a method due to Erhard Neher of constructing tame extended affine Lie algebras \cite{Neh04}, \cite{Neh11}. This construction assigns a Lie algebra $E(\fg,D,\kappa)$ to a centerless Lie torus $\fg$, where $D$ is a permissible subalgebra of skew centroidal derivations of $\fg$, and $\kappa$ is an affine cocycle on $D$; see Section \ref{tame6} for terminology. Concerning the ingredient $\kappa$, there is limited understanding of affine $2$-cocycles, and nearly all examples of EALAs are found with $\kappa=0$. 

Section \ref{EALA1} introduces the TKK construction and explains how it produces the desired  centerless Lie tori out of Jordan algebras $\bbbk_\bq^+$, $H(\bbbk_\be,-)$ and $\jj_S$. The section concludes with an analysis of an $A_1$-type elliptic extended affine Lie algebra. 
A significant challenge in constructing a Chevalley basis lies in identifying appropriate ``integral bases'' for non-zero root spaces. In this context, the complicated part of the root space associated with a non-zero element $\sg$ takes the form
$$(\star):=\bbbk L_{x^\sg}+\sum_{\tau\in\Lam}\bbbk [L_{x^{\sg+\tau}},L_{x^{-\tau}}],$$ see (\ref{dim}). Here $\Lam$ is the lattice of isotropic roots, and $L_a$ represents the left multiplication operator based on an element $a$ of the corresponding Jordan torus. The main difficulty lies in constructing an integral basis for the subspace $(\star)$, especially when working with the Jordan algebra $\jj_S$. {This will be fully resolved in Sections \ref{jordan1} and \ref{isotropic}, as we explain below.}

Section \ref{jordan1} is devoted to the Jordan torus $\jj_S$ based on a semilattice $S$ within a free abelian group $\Lam$ of finite rank. 
We begin by reviewing the multiplication operation on $\jj_S$, which relies on the involved cosets of $S$. A central element of our analysis is a function $\Gamma$, defined on these cosets, 
which helps characterize the multiplication on $\jj_S$. To facilitate our computations within the space $(\star)$, 
we establish several combinatorial results (see Lemmas \ref{final1}–\ref{coia1}).  While our primary focus is on the elliptic case, these results are proved for arbitrary nullity, making them valuable for general study of EALAs of type $A_1$.

Sections \ref{isotropic}, computes the dimension of isotropic root spaces of elliptic extended affine Lie algebras $E(\fg,D,\kappa)$ where $\fg=\hbox{TKK}(\jj)$. Here $\jj$ represents one of the rank $2$ Jordan tori  $\bbbk^+_\bq$, $H(\bbbk_\be,-)$ or $\jj_S$, $D=\hbox{SCDer}(\fg)$ is the one-dimensional algebra of skew centroidal derivations of $\fg$, and $\kappa$ is an affine $2$-cocycle. We show that if $\sg$ is a non-zero isotropic root, then the dimension of the root space $E_\sg$ is determined as follows:
$$\left\{\begin{array}{ll}
	2+\dim D^\sg+\dim C^{\sg}&\hbox{if }\jj=J_\Lam\hbox{ and }\sg\not\in2\Lam,\hbox{ or }\jj=\bbbk_\bq^+\hbox{ and }\sg\not\in\rad(f),\\
	1+\dim D^\sg+\dim C^{\sg}&\hbox{otherwise}.
\end{array}\right.
$$ 
Here $\rad(f):=\{\lam\in\Lam\mid f(\lam,\mu)=1\hbox{ for all }\mu\in\Lam\}$,	where $f(\lam,\mu)=\eta(\lam,\mu)\eta(\mu,\lam)^{-1}$ and $\eta:\Lam\times\Lam\rightarrow\bbbk^\star$ defines the product in the quantum torus $\bbbk_\bq$, namely $x^\lam\cdot x^\mu=\eta(\lam,\mu)x^{\lam+\mu}.$  Note that since $\scd(\fg)$ is one-dimensional, both $D^\sg$ and $C^\sg$ are at most one-dimensional.
It is worth mentioning that there have been several attempts to establish upper and lower bounds for the dimensions of isotropic root spaces, particulary for cases where the rank is greater than one. These attemps have often focused  on the level of the core; see, for example, \cite[Theorem 3.37]{BGKN95}, \cite[Proposition 3.31]{BGK96} and \cite[Proposition 1.4.2]{ABFP09}.

The paper concludes in Section \ref{finalsec}, which provides an explicit Chevalley basis for the elliptic extended affine Lie algebra $E(\fg,D,\kappa)$, where $D=\scd(\fg)$ and $\kappa=0$.

\markboth{S. Azam}{Chevalley Bases}

\section{\bf Preliminaries}\setcounter{equation}{0}\label{preliminaries}

All algebras are assumed to be over field $\bbbk$ of complex numbers. For a real vector space ${\mathcal V}$ equipped with a symmetric bilinear form $\fm$, and a subset
$T$ of ${\mathcal V}$, se set
$$T^\times:=\{\a\in {\mathcal V}\mid (\a,\a)\not=0\}\andd T^0:={\mathcal V}\setminus T^\times.$$ 
For $\a\in {\mathcal V}^\times$, we set 
$\a^\vee:=2\a/(\a,\a).$
For a sebset $S$ of a real vector space, we denote by $\la S\ra$ the additive subgroup generated by $S$.  

\subsection{Extended affine root systems}
We recall the definition of an extended affine root system from \cite{AABGP97}.
\begin{DEF}\label{azam1}
	An {\it extended affine root system} is a triple
	$({\mathcal V},\fm,R)$ where ${\mathcal V}$ is a finite dimensional real vector space, $\fm$ is a symmetric positive semi-definite bilinear form, and $R$ is a subset of ${\mathcal V}$ satisfying the following seven axioms:
	
	(R1) $0\in R$,
	
	(R2) $R=-R$,
	
	(R3) $R$ spans ${\mathcal V}$,
	
	(R4) $\a\in R^\times\Rightarrow 2\a\not\in R$,
	
	(R5) $R$ is discrete in ${\mathcal V}$,
	
	(R6) the {\it root string property}: for $\a\in R^\times$ and $\beta\in R$, there exist non-negative integers $u,d$ such that
	$$\{i\in\bbbz\mid\beta+i\a\in R\}=\{-d,\ldots,u\}
	$$
	with $(\beta,\a^\vee)=d-u$,
	
(R7) elements of $R^0$ are {\it non-isolated}, i.e., for $\sg\in R^0$ there exists $\a\in R^\times$ with $\a+\sg\in R$.

	(R8)  $R^\times$ is {\it connected}, meaning that $R^\times$ cannot be written in the form $R^\times=R_1\cup R_2$ with
$(R_1,R_2)=\{0\}$ where $R_1\not=\emptyset$ and $R_2\not=\emptyset$. 
\end{DEF}

\begin{parag}\label{pt1}
Let $(R, \fm, {\mathcal V})$ be an extended affine root system. Define $\bar{\mathcal V} := {\mathcal V}/ {\mathcal V}^0$, with $\bar{\;} : {\mathcal V} \rightarrow \bar{{\mathcal V}}$ as the canonical map. 
Then, $\bar{R}$, the image of $R$ under $\bar{\;}$, is an irreducible finite root system in $\bar{{\mathcal V}}$, with respect to the form induced by $\fm$ on $\bar{{\mathcal V}}$. The {\it type} and {\it rank} of $R$ are defined as those of $\bar{R}$, and the {\it nullity} of $R$ is defined as the dimension of ${\mathcal V}^0$. Throughout this work, $R$ is {\it always of reduced type}.

{It follows that $R$ contains a finite root system $\dot{R}$ that is isomorphic to $\bar{R}$ under $\bar{\;}$. Let $\rds$ and $\rdl$ be the sets of short and long roots of $\rd$, respectively.  Then the root system $R$ can be expressed as
\[
R = (S + S) \cup (\rds + S) \cup (\rdl + L) \quad \text{with} \quad R^0 = S + S,
\]
for certain subsets $S$ and $L$ of ${\mathcal V}^0$, referred to as {\it semilattices}. We recall that a semilattice $S$ in ${\mathcal V}^0$ is a discrete spanning subset of ${\mathcal V}^0$ satisfying $0\in S$, and $S\pm 2\la S\ra\sub S$. If $\dot{R}$ is simply laced, we treat all roots as short, meaning $\dot{R}^\times = \rds$, and interpret $\rdl + L$ as an empty set. The semilattices $S$ and $L$ interact in such a way that
$$S+\la L\ra\sub S,\andd L+k\la S\ra\sub L,$$ where $k$ is the maximum number of multiple edges between two nodes appearing in the Dynkin diagram of $\rd$.
When $\dot{R}$ is of type $X$, we denote $R$ as
$
R = R(X, S)$ or $R = R(X, S, L),
$ depending on whether $R$ is simply laced or not. For details about semilattices, we refer the reader to \cite[Chapter II]{AABGP97}
and \cite{Az97}.}
\end{parag}

\subsection{\bf Extended affine Lie algebras}\setcounter{equation}{0}
We briefly recall the definition and some properties of extended affine Lie algebras.
\begin{DEF}\label{ea1}
	An {\it extended affine Lie algebra} is  a triple $(E,\fm,\hh)$ where $E$ is a Lie algebra, $\hh$ is a subalgebra of $E$ and $\fm$ is a bilinear form on $E$  satisfying the following five axioms:
	
	(EA1) The form $\fm$ on $E$ is symmetric, non-degenerate and invariant.
	
	(EA2)  $\hh$ is a finite dimensional splitting Cartan subalgebra of $E$. This means that $E=\sum_{\a\in\hh^\star}E_\a$ where $E_\a=\{x\in E\mid [h,x]=\a(h)x\hbox{ for all }h\in\hh\}$ and $E_0=\hh$.

	Let $R$ be the set of roots of $E$, that is, $R=\{\a\in\hh^\star\mid E_\a\neq\{0\}\}$.
	Axioms (EA1)-(EA2) imply that the restriction of $\fm$ to $\hh$ is non-degenerate, allowing it to transferred
	to $\hh^\star$ via $(\a,\beta):=(t_\a,t_\beta)$ where $t_\a\in\hh$
	is the unique element satisfying $\a(h)=(h,t_\a)$, $h\in\hh$.
	Let 
	$$R^0=\{\a\in R\mid (\a,\a)=0\}\andd R^\times=R\setminus R^0.$$
	Then $R=R^0\uplus R^\times$ is regarded as the decomposition of roots into
	{\it isotropic} and {\it non-isotropic} roots, respectively. Let $V:=\hbox{Span}_\bbbr R$
	and ${\mathcal V}^0:=\hbox{Span}_{\bbbr}R^0$.
	
	(EA3) For $\a\in R^\times$, the adjoint action $\ad x$ is locally nilpotent for $x\in E_\a$.
	
	(EA4) The $\bbbz$-span of $R$ is a full lattice in ${\mathcal V}$, meaning that ${\mathcal V}$ has a $\bbbr$-basis which is also a basis for the $\bbbz$-span of $R$.
		
	(EA5) $R$ is irreducible, that is, $R^\times$ is connected and elements of $R^0$ are non-isolated.
	\end{DEF}
	
Here we record some basic facts about the extended affine Lie algebra $(E,\fm,\hh)$ relevant to this work. For a detailed study of extended affine Lie algebras and their root systems, refer to \cite{AABGP97} and \cite{Neh11}. Let $R$ be the root system of $E$. It turns out that $R$ is an extended affine root system in the sense of definition \ref{azam1}. The {\it type}, {\it rank}, and {\it nullity} of $E$ are defined as the same for $R$. As we have already stated, $R$ and so $E$ is always assumed to be reduced.

From \cite[Chapter I.\S 1]{AABGP97}, we have 
\begin{equation}\label{khat0}
[E_\a,E_{-\a}]=\bbbk t_\a,\quad(\a\in R),
\end{equation}
\begin{equation}\label{khat1}
	(E_\a,E_\beta)=\{0\}\hbox{ unless }\a+\beta=0,\quad(\a,\beta\in R).
	\end{equation} 
	Also from \cite[Remark 1.5]{Az06}, we have
	\begin{equation}\label{khat4}
		[E_\a, E_\beta]\not=\{0\},\qquad(\a\in R^\times,\beta,\a+\beta\in R).
		\end{equation}

The {\it core} $E_c$ of $E$ {is by definition the subalgebra} generated by non-isotropic root spaces of $E$. It follows that $C_E(E_c)=E_c^\perp$, where $C_E(E_c)$ is the centralizer of the core in $E$, and $E_c^\perp$ is the orthogonal complement of $E_c$ in $E$ with respect to the form. The extended affine Lie algebra $E$ is called {\it tame} if $E_c^\perp=Z(E_c)$.

\begin{lem}\label{naz4}
For $\sg\in R^0$ and $\a, \a+\sg\in R^\times$, we have
 $([x_{\a+\sg},x_{-\a}],[x_{-\a-\sg},x_{\a}])\not=0.$
 \end{lem}
 
 \proof From the Jacobi identity, the fact that $-2\a-\sg$ is not a root, and invariance of the form, we have 
\begin{eqnarray*}
 ([x_{\a+\sg},x_{-\a}],[x_{-\a-\sg},x_{\a}])
 &=& (x_{\a+\sg},[x_{-\a},[x_{-\a-\sg},x_{\a}]])\\
 &=&
 -(x_{\a+\sg}, [x_{-\a-\sg},[x_\a,x_{-\a}]])\\
 &=&
 -(x_{\a+\sg},[x_{-\a-\sg},t_\a])\\
 &=&
 -([x_{\a+\sg},x_{-\a-\sg}],t_\a)\\
 &=&-(t_{\a+\sg},t_\a)=-(\a,\a)\not=0.
 \end{eqnarray*}
 \qed

\section{construction of extended affine Lie algebras}\setcounter{equation}{0}\label{tame6}
In this section, we review a construction of tame extended affine Lie algebras due to Erhard Neher, \cite{Neh11}. We begin with the definition of a Lie torus.

\subsection{Lie tori}\label{Lie tori}
Assume that $\Delta$ is an irreducible finite root system and $\Lam$ is a free abelian group of finite rank.
\begin{DEF}\label{def1}
	A Lie $\Lam$-torus of type $\Delta$ is a Lie algebra $\fg$ over $\bbbk$ that satisfies the following conditions (LT1)-(LT4):\\
	(\textbf{LT1}) $\fg$ has  a compatible $(Q\times\Lam)$-grading $\fg=\bigoplus_{(\a,\lam)\in Q\times\Lam} \fg_\a^\lam$m with $\fg_\a^\lam=0$ if $\a\notin\Delta,$
	and $[\fg_\a^\lam,\fg_\beta^\mu]\subset \fg_{\a+\beta}^{\lam+\mu}.$
	
	\noindent (\textbf{LT2}) For $\a\in\Delta^\times$ and $\lam\in\Lam$ we have
	$\dim \fg_\a^\lam\leq1$, $\dim \fg_\a^0=1$ if $\a\in\Delta_{\ind}$. If
	 $\dim \fg_\a^\lam=1$ then there exist elements $e_\a^\lam\in \fg_\a^\lam$ and $f_\a^\lam\in \fg_{-\a}^{-\lam}$ such that
		$$\fg_\a^\lam=\bbbk e_\a^\lam,\quad\fg_{-\a}^{-\lam}=\bbbk f_\a^\lam,
		\andd 
			[[e_\a^\lam,f_\a^\lam],x^\mu_\beta]=\langle\beta,\a^\vee\rangle x_\beta^\mu$$
			for $\beta\in\Delta,$ $\mu\in\Lam$, $x_\beta^\mu\in \fg_\beta^\mu.$

	\noindent (\textbf{LT3}) For $\lam\in\Lam$, $\fg_0^\lam=\sum_{\a\in\Delta^\times,\mu\in\Lam} [\fg_\a^\mu,\fg_{-\a}^{\lam-\mu}]$.\\
	(\textbf{LT4}) $\Lam=\langle\supp_\Lam(\fg)\rangle$, where $\supp_\Lam(\fg)=\{\lam\in\Lam\mid \fg_\a^\lam\neq 0\;\text{for some} \;\a\in\Delta\}$.
	
The rank and the \textit{nullity} of $\fg$ are defined as the rank of $\Delta$ and the nullity of $\Lam$, respectively. A  Lie torus $\fg$ is called {\it centreless} if $\fg$ has trivial center.		
\end{DEF}

\subsection{A construction of extended affine Lie algebras}\label{const}
Assume that $\fg$ is a centerless $\Lam$-Lie torus.
Recall that the \textit{centroid} of $\fg$, denoted  $\cent_\bbbk(\fg)$, is the set of linear endomorphisms that commute with left and right multiplications by elements of $\fg$, i.e.,
$$\cent_\bbbk(\fg)=\{\chi\in\mathrm{End}(\fg)\mid\chi([x,y])=[\chi(x),y]=[x,\chi(y)],\forall x,y\in \fg\}.$$
One knows from  \cite{Neh04} that 
\begin{equation}\label{ifo1}
	\cent_\bbbk(\fg)=\bigoplus_{\mu\in\Gamma}\bbbk\chi^\mu,
	\end{equation}
where $\Gamma$ is the subgroup of $\Lam$ consisting of $\lam\in\Lam$ for which $\cent_\bbbk(\fg)^\lam\neq 0$. The group $\Gamma$ is called the {\it central grading group} of $\fg$. $\chi^\mu$ acts on $\fg$ as an endomorphism of degree $\mu$ with the property that $\chi^\mu\chi^\nu=\chi^{\mu+\nu}$.

For $\theta\in\text{Hom}_\bbbz(\Lam,\bbbk)$, the {\it degree derivation} $\partial_\theta$ of $\fg$ is defined by
$$\partial_\theta(x^\lam)=\theta(\lam)x^\lam\;\;\text{for}\;\lam\in\Lam,x^\lam\in \fg^\lam.$$
We denote by  $\mathcal{D}$ the set of all degree derivations, and by $\mathrm{CDer}_\bbbk(\fg)$
the set of \textit{centroidal derivations} of $\fg$, namely
$$\mathrm{CDer}_\bbbk(\fg):=\cent_\bbbk(\fg)\mathcal{D}=\bigoplus_{\mu\in\Gamma}\chi^\mu\mathcal{D},$$
which is a $\Gamma$-graded subalgebra of the derivation algebra
$\mathrm{Der}_\bbbk(\fg)$ of $\fg$ with
\begin{equation}\label{eq7}
	[\chi^\mu\partial_\theta,\chi^\nu\partial_\psi]=\chi^{\mu+\nu}(\theta(\nu)\partial_\psi-\psi(\mu)\partial_\theta).
\end{equation}

Next, the set
\begin{eqnarray*}
	\mathrm{SCDer}_\bbbk(\fg)&:=&\{d\in\mathrm{CDer}_\bbbk(\fg)\mid(d(x),x)_\fg=0\;\text{for all}\;x\in \fg\}\\
	&=&\bigoplus_{\mu\in\Gamma}\mathrm{SCDer}_\bbbk(\fg)^\mu=\bigoplus_{\mu\in\Gamma}\chi^\mu\{\partial_\theta\in\dd\mid\theta(\mu)=0\},
\end{eqnarray*}
is a $\Gamma$-graded subalgebra of $\mathrm{CDer}_\bbbk(\fg)$, called the algebra of {\it skew centroidal derivations} of $\fg$. Note that $\mathrm{SCDer}_\bbbk(\fg)^0=\mathcal{D}.$

For a graded subalgebra $D=\sum_{\mu\in\Gamma}D^\mu$ of
$\mathrm{SCDer}_\bbbk(\fg)$ denote its graded dual by
$ D^{gr^\star}=\sum_{\mu\in\Gamma}(D^\mu)^\star$ with grading
$(D^{gr^\star})^\mu=(D^{-\mu})^\star$, and consider it as a $D$-module by
$$(d.\varphi)(d^{'})=\varphi([d^{'},d])\:\text{for}\;d^{'},d\in D,\varphi\in D^{gr*},$$
where $\varphi\in(D^\mu)^\star$ is viewed as a linear map on $D$ by $\varphi|_{D^\nu}=0$ for $\nu\neq\mu$. 

To construct an extended affine Lie algebra we introduce two more ingredients. 
A $\Gamma$-graded subalgebra $D=\bigoplus_{\mu\in\Gamma}D^\mu$ of  $\mathrm{SCDer}_\bbbk(\fg)$ is called
{\it permissible} if the canonical evaluation map $\text{ev}:\Lam\rightarrow(D^0)^\star$ defined by
$\text{ev}(\lam)(\partial_\theta)=\theta(\lam)$, $\lam\in\Lam,$
is injective {and has discrete image}.  Lastly, let $\kappa$ be a bilinear map satisfying
$$\begin{array}{c}
	\kappa(d,d)=0,\;\;{\sum_{(i,j,k)\circlearrowleft}\kappa([d_i,d_j],d_k)=\sum_{(i,j,k)\circlearrowleft}d_i\cdot\kappa (d_j,d_k)},\vspace{2mm}\\
	\kappa(D^{\mu_1},D^{\mu_2})\subseteq(D^{-\mu_1-\mu_2})^\star\;\;\text{and}\;\;\kappa(d_1,d_2)(d_3)=\kappa(d_2,d_3)(d_1),\vspace{2mm}\\
	\kappa(D^0,D)=0,
\end{array}
$$
for $d,d_1,d_2,d_3\in\ D$.
{Note that by $(i,j,k)\circlearrowleft$, we mean $(i,j,k)$ is a cyclic permutation of $(1,2,3)$}. $\kappa$ is called an {\it affine cocycle} on $D$.

Assume now that $\fg$ is a centerless Lie torus, $D$ is a permissible subalgebra of
$\scd(\fg)$ and $\kappa$ is an affine cocycle on $D$. Set
\begin{equation}\label{rok1}
	E=E(\fg,D,\kappa):=\fg\oplus D^{gr*}\oplus D.
	\end{equation}
Then $E$ is a Lie algebra with the bracket
\begin{eqnarray*}
	[x_1+c_1+d_1,x_2+c_2+d_2]&=&([x_1,x_2]_\fg+d_1(x_2)-d_2(x_1))\\
	&+&(c_D(x_1,x_2)+d_1.c_2-d_2.c_1+\kappa(d_1,d_2))\\
	&+&[d_1,d_2]
\end{eqnarray*}
for $x_1,x_2\in \fg,c_1,c_2\in D^{gr*},d_1,d_2\in D$. Here $[\:,\:]_\fg$ denotes the Lie bracket on $\fg$, $[d_1,d_2]=d_1d_2-d_2d_1$, and $c_D:\fg\times\fg\rightarrow D^{gr*}$ is defined by
$$c_D(x,y)(d)=(d(x)|y)\;\text{for all}\;x,y\in \fg,d\in D.$$
Next, define a bilinear form $\fm$ on $E$ by
\begin{equation}\label{ira1}
	(x_1+c_1+d_1,x_2+c_2+d_2)=(x_1,x_2)_\ll+c_1(d_2)+c_2(d_1).
\end{equation}
It follows that $\fm$ is symmetric, invariant and non-degenerate. In fact the Lie algebra $(E,\fm,\hh)$ constructed from the data $\fg,$ $D$ and $\kappa$ is a tame extended affine Lie algebra, where 
$$\hh=\fg^0_0\oplus(D^0)^\star\oplus D^0,$$
and conversely any tame extended affine Lie algebra arises this way, see \cite[Theorem 16]{Neh04}. To indicate the dependence of $E$ on $\fg$, $D$ and $\kappa$, we write
$E=E(\fg,D,\kappa)$. Note that
$$E_c=\fg\oplus {D^{gr}}^\star\andd
	Z(E_c)={D^{gr}}^\star.
$$

\begin{rem}\label{tame7}
	(i)	From $	Z(E_c)={D^{gr}}^\star$, we have
	\begin{equation}\label{tame8}
		Z(E_c)\sub\hh\Longleftrightarrow D=D^0\Longleftrightarrow Z(E_c)={D^0}^\star.
	\end{equation}	
	
	(ii) Almost all examples of extended affine Lie algebras found in the literature have typically $D=D^0$ and $\kappa=0$, as one can see for example in \cite[Chapter III]{AABGP97}, \cite{BGKN95}, \cite{BGK96} and \cite{H-KT90}.
\end{rem}

\section{Extended affine Lie algebras of type $A_1$}\setcounter{equation}{0}\label{EALA1}
 We begin by recalling the TKK-construction which associates to any Jordan algebra a Lie algebra.
 \subsection{TKK construction}\label{tkk}

Let $\jj$ be a Jordan algebra, and set $\mathrm{Inder}(\jj):=\{L_x\mid x\in \jj\}$ where $L_x$ is the operator on $\jj$ defined by $L_xy=xy$, for $y\in\jj$. Next set
$$\mathrm{Instrl}(\jj):=\sum_{x\in \jj}L_x+\{\sum_i[L_{x_i},L_{y_i}]\mid x_i,y_i\in\jj\},$$ 
and
$$ \mathrm{TKK}(\jj):=\jj\oplus\mathrm{Instrl}(\jj)\oplus\bar{\jj},$$
where $\bar\jj$ is a copy of $\jj$.
Then $\mathrm{Instrl}(\jj)$ and $\mathrm{TKK}(\jj)$ are Lie algebras with the brackets given by
$$[L_x+C,L_y+D]=[L_x,L_y]+L_{Cy}-L_{Dx}+[C,D],$$
for $x,y\in\jj,C,D\in\mathrm{Inder}(\jj)$, and
\begin{equation}\label{naz1}
\begin{array}{c}
	[x_1+\bar{y}_1+E_1,x_2+\bar{y}_2+E_2]=-E_2x_1+E_1x_2-\overline{\bar E_2 y}+\overline{\bar E_1y_2}\\
	\qquad\qquad\qquad\qquad\qquad\qquad\qquad+x_1\triangle y_2-x_2\triangle y_1+[E_1,E_2],
	\end{array}
\end{equation}
for $x_i\in\jj,\bar y_i\in\bar{\jj}$, and $E_i\in\mathrm{Instrl}(\jj)$, where $x\triangle y=L_{xy}+[L_x,L_y]$.
Here $\bar{}:\mathrm{Instrl}(\jj)\rightarrow\mathrm{Instrl}(\jj)$ is an involution define by 
 $\overline{L_x+D}=-L_x+D$. The Lie algebra $\mathrm{TKK}(\jj)$ is called TKK {\it Lie algebra} of $\jj$.

Next, let $\Lam$ be a free abelian group of rank $\nu$. The Jordan algebra $\jj$ is called a $\Lam$-Jordan torus if $\jj=\bigoplus_{\lam\in\Lam}\jj^\lam$ with
$\dim \jj^\lam\leq 1$ for each $\lam$, and $\Lam$ is generated by $S:=\{\lam\in\Lam\mid \jj^\lam\not=\{0\}\}$.
Set
$$\fg:=\mathrm{TKK}(\jj)$$
Then $\fg$ is $\Lam$-graded with
\begin{equation}\label{eq9}
\fg^\lam=	\mathrm{TKK}(\jj)^\lam=\jj^\lam\oplus\mathrm{Instrl}(\jj)^\lam\oplus\overline{\jj^\lam},
\end{equation}
for $\lam\in\Lam$, where $\mathrm{Instrl}(\jj)^\lam=L_{\jj^\lam}\oplus\sum_{\mu+\nu=\lam}[L_{\jj^\mu},L_{\jj^\nu}]$.

The Lie algebra $\fg$ can be equipped with a symmetric invariant non-degenerate form as follows. We fix a basis $\{x^\lam\in\jj^\lam\mid \lam\in S\}$ for $\jj$. For the sake of notation if $\lam\in\Lam\setminus S$, we write $x^\lam$ and interpret it as zero.
Consider the linear map $\ep:\jj\rightarrow\bbbk$ induced by $\ep(1)=1$ and $\ep(x^\lam)=0$ if $\lam\not=0$. Then
$(x,y):=\ep(xy)$ defines a symmetric invariant non-degenerate form on $\jj$. This then defines a form on $\mathrm{Instrl}(\jj)$ by 
$(D,[L_x,L_y])=(Dx,y)$ for $D\in\mathrm{Instrl(\jj)}$, $x,y\in\jj$. Finally, this form extends to a symmetric invariant non-degenerate form on $\fg$ by
{\small$$(x_1+L_{y_1}+D_1+\bar{z}_1, x_2+L_{y_2}+D_2+\bar{z}_2):=(x_1,z_2)+(x_2,z_1)+(y_1,y_2)+(D_1,D_2),$$}
for $x_i,y_i,z_i\in\jj$, $D_i\in\mathrm{Instrl}(\jj)$.

\subsection{Elliptic extended affine Lie algebras of type $A_1$}\setcounter{equation}{0}\label{elliptic1}

We recall from \cite[Corollary 7.9]{Yos00} that any centerless Lie torus of type $A_1 $ takes the form
$\fg=\mathrm{TKK}(\jj)$ where $\jj$ represents one of the following Jordan algebras:
	$\bbbk_\bq^+$, $H(\bbbk_\be,-)$, $\jj_S$, or $\mathbb{A}_t$. {Since this work focuses on the elliptic case (nullity 2), we restrict our discussion to the first three Jordan algebras. The Jordan algebra
$\mathbb{A}_t$ 
does not arise in the classification of $A_1$-type extended affine Lie algebras of nullity $1$ or $2$; it only appears for nullity $\geq 3$, which falls outside the scope of the present study.}

Consider the centerless Lie torus $\fg=\mathrm{TKK}(\jj)$ where $\jj$ is one the Jordan algebras (Jordan tori) mentioned above. Let $\Lam=\bbbz^\nu.$
Set $\dot\hh=\bbbk L_1$ and define $\da\in{\dot\hh}^\star$ by $\da(L_1)=1$. 
Set $Q=\Z\da$.
Equip $\fg$ with the $Q$-grading
\begin{equation}\label{eq2}
	\fg_\beta=\left\{\begin{array}{ll}
		\jj&\hbox{if }\beta=\da,\\
		\bar{\jj}&\hbox{if }\beta=-\da,\\
		\mathrm{Instrl}(\jj)&\hbox{if }\beta=0,\\
		\{0\}&\hbox{otherwise.}
	\end{array}\right.
\end{equation}

Next, let $D$ be a permissible subalgebra of the skew centroidal derivations $\scd(\fg)$  of $\fg$ and $\kappa$ be an affine cocycle on $D$. Let $C$ be the graded dual of $D$.
Set $E=E(\fg,D,\kappa)$  as in Section \ref{tame6}. Then
$$E:=\fg\oplus C\oplus D\andd\hh:=\bbbk L_1\oplus C^0\oplus D^0.$$
By \cite[Theorem 16]{Neh04}, $(E,\fm,\hh)$ is a tame extended affine Lie algebra of type $A_1$ with root system $R\sub\Lam\cup(\pm\da+\Lam)$, and any tame extended affine Lie algebra of type $A_1$ is up to isomorphism of this form. Moreover, for $\beta\in R,$
\begin{equation}\label{dim}
	E_\beta=\left\{
	\begin{array}{ll}
		\bbbk L_1\oplus D^0\oplus C^0&\hbox{if }\beta=0,\\
		\bbbk x^\sg&\hbox{if }\beta=\da+\sg,\;\sg\not=0,\\
		\bbbk\bar{x}^\sg&\hbox{if }\beta=-\da+\sg,\;\sg\not=0,\\
		\bbbk L_{x^\sg}+\sum_{\tau\in\Lam}\bbbk [L_{x^{\sg+\tau}},L_{x^{-\tau}}]\oplus D^\sg\oplus C^\sg&\hbox{if }\beta=\sg,\;\sg\not=0,\\
	\end{array}
	\right.
\end{equation}

We note that since $(L_1,L_1)=\ep(1)=1$, it follows that $(\dot\a,\dot\a)=1$.

\begin{rem}\label{rok5}
(i) From the way the bilinear form $\fg$ is defined on $\fg=\hbox{TKK}(\jj)$ and on $E=E(\fg,D,\kappa)$, we have
$(L_1,L_1)=(1,1)=\ep(1)=1$ and so $(\dot\a,\dot\a)=1$. 

(ii) If $\rank\;\Lam=2$, then $\dim\Hom_\bbbz(\Lam,\bbbk)=2$ and so $\dim\scd(\fg)^\sg\leq 1$ for each $0\not=\sg\in\Lam$.
Therefore, $\dim D^\sg,\dim C^\sg\leq 1$ for $0\not=\sg\in\Lam$.
\end{rem}

\section{The Jordan torus $\jj_S$}\setcounter{equation}{0}\label{jordan1}
In this section, we examine the $\Lam$-Jordan torus $\jj_S$ based on $S$ where $S$ is a semilattice in a free abelian group $\Lam$ of finite rank. The results presented here are not limited to nullity $2$; they apply to any arbitrary nullity.

\subsection{The function $\Gamma$}
Let $S$ be a semilattice in $\Lam=\bbbz^\nu$, that is, $S=\cup_{i=0}^mS_i$ where
$S_i=\tau_i+2\Lam$ and $\tau_1,\ldots,\tau_m$ represent distinct cosets of $2\Lam$ in $\Lam$ with $\tau_0=0$. 
For $\sg\in S$ consider the symbol $x^\sg$ and set $\jj=\jj_S=\sum_{\sg\in S}\bbbk x^\sg$. By convention, we write $\jj=\sum_{\sg\in\Lam}\bbbk x^\sg$, where we interpret $x^\sg=0$ if $\sg\not\in S$. Then $\jj$ is a $\Lam$-Jordan torus with the multiplication
\begin{equation}\label{rokn2}
	x^\sg\cdot x^\tau=\left\{\begin{array}{ll}
	x^{\sg+\tau}&\hbox{if }\sg,\tau\in S_0\cup S_i,\; 0\leq i\leq m,\\
	0&\hbox{otherwise}.
\end{array}\right.
\end{equation}

We choose cosets $S_i=\tau_i+2\Lam$, $i=m+1,\ldots, 2^\nu-1$ of $2\Lam$ in $\Lam$ such that $\Lam=\cup_{i=0}^{2^\nu-1}S_i$. Now, we define a symmetric function $\Gamma:\ss\times\ss\rightarrow\{0,1\}$ where $\ss=\{S_i\mid 0\leq i\leq 2^\nu-1\}$ by
\begin{equation}\label{tem1}
	\Gamma(S_i,S_j)=\left\{\begin{array}{ll}
	0&\hbox{if $i$ or $j\geq m+1$}\\
	0&\hbox{if }i\not=j >0\\
	1&\hbox{otherwise.}
	\end{array}\right.
	\end{equation}
In particular, for $0\leq i\not=j\leq m$, we have
$$\Gamma(S_0,S_0)=\Gamma(S_0,S_i)=\Gamma(S_i,S_i)=1\andd \Gamma(S_i,S_j)=0.$$
This gives 
\begin{equation}\label{tem2}
\Gamma (S_i,S_j)=0,\qquad (S_i+S_j\not\in\ss,\; 0\leq i,j\leq m),
\end{equation}
since in this case $i,j$ can not be equal and none can be zero.
By convention we set
$$\Gamma(\sg,\tau)=\Gamma(S_i,S_j)=\Gamma(\sg,S_j)\qquad(\sg\in S_i,\tau\in S_j).$$ 
Then the map $\Gamma$ can be describe the multiplication (\ref{rokn2}) as 
$$ \begin{array}{c}
x^\sg\cdot x^\tau=\Gamma(\sg,\tau)x^{\sg+\tau}=\Gamma(S_i,S_j)x^{\sg+\tau},\vspace{2mm}\\
(\sg\in S_i,\tau\in S_j,\;0\leq i,j\leq 2^\nu-1).
\end{array}
$$
We note that since for each $i$, $S_0+S_i=S_i$ and $S_i+S_i=S_0$, we have for $0\leq i,j\leq m$,
\begin{equation}\label{rokn3}
	\Gamma(S_0, S_0+S_i)=\Gamma(S_i,S_j+S_j)=\Gamma(S_i,S_0+S_i)=1.
	\end{equation}

\begin{lem}\label{rokn1}
	For $\lam_i\in S_i$ we have $x^{\lam_i}\cdot(x^{\lam_j}\cdot x^{\lam_k})=\Gamma(S_j,S_k)\Gamma(S_i,S_j+S_k)x^{\lam_i+\lam_j+\lam_k}.$
	\end{lem}
\proof
If either of $\lam_i,\lam_j,\lam_k$ is in $\Lam\setminus S$, then $x^{\lam_i}\cdot(x^{\lam_j}\cdot x^{\lam_k})=0$ and 
$\Gamma(S_j,S_k)\Gamma(S_i,S_j+S_k)=0$, so we are done in this case.
Assume now that 
 $\lam_i\in S_i$, $\lam_j\in S_j$ and $\lam_k\in S_k$, where $0\leq i,j,k\leq m$. Then
 $$x^{\lam_i}\cdot(x^{\lam_j}\cdot x^{\lam_k})=\Gamma(S_j,S_k)x^{\lam_i}\cdot x^{\lam_j+\lam_k}=
 \Gamma(S_j,S_k)\Gamma(S_i,S_j+S_k)x^{\lam_i+\lam_j+\lam_k}.$$\qed

\subsection{Lie Brackets of left operators}
We establish several results on brackets containing left operators which facilitate in computing dimensions of isotropic root spaces.
\begin{lem}\label{final1}
	(i) $[L_{x^{\sg}},L_{x^{\tau}}](x^\gamma)=0$ if at least one of $\sg,\tau$ or $\gamma$ belongs to $S_0$.
	
	(ii) $[L_{x^{\sg}},L_{x^{\tau}}]=0$ if $\sg,\tau\in S_i$  for some $i$.
\end{lem}

\proof
(i) Let $\sg\in S_0$. We may assume that $\tau,\gamma\in S$. If $\tau+\gamma\not\in S$, then
$\tau+\gamma+\sg\not\in S$ and by (\ref{tem1}) and (\ref{tem2}),
	$$[L_{x^{\sg}},L_{x^{\tau}}](x^\gamma)=
\big(\Gamma(\tau,\gamma)\stackrel{=0}{\overbrace{\Gamma(\sg,\tau+\gamma)}}-\Gamma(\sg,\gamma)\stackrel{=0}
{\overbrace{\Gamma(\tau,\sg+\gamma)}}\big){x^{\sg+\tau+\gamma}}=0.$$
If $\tau+\gamma\in S$, then
\begin{eqnarray*}
	[L_{x^{\sg}},L_{x^{\tau}}](x^\gamma)
		&=&\big(\Gamma(\tau,\gamma)\Gamma(S_0,\tau+\gamma)-\Gamma(S_0,\gamma)\Gamma(\tau,S_0+\gamma)\big)x^{\sg+\tau+\gamma}\\
	&=&\big(\Gamma(\tau,\gamma)-\Gamma(\tau,S_0+\gamma)\big)x^{\sg+\tau+\gamma}\\
	&=&\big(\Gamma(\tau,\gamma)-\Gamma(\tau,\gamma)\big)x^{\sg+\tau+\gamma}=0.
	\end{eqnarray*}
		The case $\tau\in S_0$ holds by symmetry.

If $\gamma\in S_0$, then
\begin{eqnarray*}
	[L_{x^{\sg}},L_{x^{\tau}}](x^\gamma)
	&=&\big(\Gamma(\tau,S_0)\Gamma(\sg,\tau+S_0)-\Gamma(\sg,S_0)\Gamma(\tau,\sg+S_0)\big)x^{\sg+\tau+\gamma}\\
	&=&(\Gamma(\sg,\tau)-\Gamma(\tau,\sg))x^{\sg+\tau+\gamma}=0.
\end{eqnarray*}

(ii) 
 Assume $\sg,\tau\in S_i$ and $\gamma\in S_j$ for some $i,j$. Then
\begin{eqnarray*}
	[L_{x^{\sg}},L_{x^{\tau}}](x^\gamma)&=&
	\big(\Gamma(S_i,S_j)\Gamma(S_i,S_i+S_j)-\Gamma(S_i,S_j)\Gamma(S_i,S_i+S_j)\big)x^{\sg+\tau+\gamma}\\
	&=&0.
\end{eqnarray*}
\qed

\begin{lem}\label{new09}
	Let $\sg,\tau,\gamma\in S$ and $\lam_0,\lam'_0,\lam''_0\in S_0$ with $\lam_0+\lam'_0+\lam''_0=0$. Then
	$[L_{x^{\sg}}, L_{x^{\tau}}](x^\gamma)=[L_{x^{\sg+\lam_0}}, L_{x^{\tau+\lam'_0}}](x^{\gamma+\lam''_0}).$ In particular,
	if $\tau\in S_j$, then
	$$[L_{x^{\sg}}, L_{x^{\tau}}]=[L_{x^{\sg+\tau+\tau_j}}, L_{x^{-\tau_j}}].$$
\end{lem}

\proof Let $\sg\in S_i$, $\tau\in S_j$ and $\gamma\in S_k$. Then
\begin{eqnarray*}
&&[L_{x^{\sg+\lam_0}}, L_{x^{\tau+\lam'_0}}](x^{\gamma+\lam''_0})\\	
&&\qquad\qquad
	=\big(\Gamma(S_j+S_0,S_k+S_0)\Gamma(S_i+S_0,S_j+S_k+S_0))\\
	&&\qquad\qquad\qquad-\Gamma(S_i+S_0,S_k+S_0)\Gamma(S_j+S_0,S_i+S_k+S_0)\big)x^{\sg+\tau+\gamma}\\
	&&\qquad\qquad=
	\big(\Gamma(S_j,S_k)\Gamma(S_i,S_j+S_k))-\Gamma(S_i,S_k)\Gamma(S_j,S_i+S_k)\big)x^{\sg+\tau+\gamma}\\
	&&\qquad\qquad=
	[L_{x^{\sg}}, L_{x^{\tau}}](x^\gamma).
	\end{eqnarray*}
This proves the first assertion.

	The second assertion  is immediate from the first one as $\tau+\tau_j\in S_0.$\qed

\begin{lem}\label{let}
	Suppose $\lam_i+\lam_j+\lam_k\in S_0$. Then we have
	$[L_{x^{\lam_i}},[L_{x^{\lam_j}},L_{x^{\lam_k}}]]=0.$
\end{lem}	

\proof Let $\lam_t\in S_t$ for $t=i,j,k$. If at least one of $\lam_i,\lam_j,\lam_k$ is in $S_0$, then using Lemma \ref{final1}(i) and Jacobi identity we are done. Suppose next that at least two of $\lam_i$, $\lam_j$, $\lam_k$ belong to the same coset. If $\lam_j$ and $\lam_k$ belong to the same coset then by Lemma \ref{final1}(ii) we are done. If $\lam_i,\lam_j$, or $\lam_i,\lam_k$ belong to the same coset, then form $\lam_i+\lam_j+\lam_k\in S_0$, we conclude that $\lam_k$ or $\lam_j$ belongs to $S_0$ and so we are done by Lemma \ref{final1}(i). 

The above discussion shows that we may assume $i,j,k$ are distinct and none is zero. Note that by assumption
$S_i+S_j+S_k=\lam_i+\lam_j+\lam_k+2\Lam=S_0$, so $S_j+S_k=S_i$. We now compute
$[L_{x^{\lam_i}},[L_{x^{\lam_j}},L_{x^{\lam_k}}]](x^{\lam_u})=J_1-J_2$ for $\lam_u\in S$, where
\begin{eqnarray*}
&&	J_1:={x^{\lam_i}}\cdot ([{x^{\lam_j}}\cdot {x^{\lam_k}}]\cdot x^{\lam_u})\\
	&&\quadd\qquad=
	\big(\Gamma(S_k,S_u)\Gamma(S_j,S_k+S_u)-\Gamma(S_j,S_u)\Gamma(S_k,S_j+S_u)\big)x^{\lam_i}\cdot x^{\lam_j+\lam_k+\lam_u}\\
	&&\quadd\qquad=\Gamma(S_i,S_j+S_k+S_u)\big(\Gamma(S_k,S_u)\Gamma(S_j,S_k+S_u)\\
	&&\quadd\qquad\quadd\qquad -\Gamma(S_j,S_u)\Gamma(S_k,S_j+S_u)\big)x^{\lam_i+\lam_j+\lam_k+\lam_u}\\
		&&\quadd\qquad=\Gamma(S_i,S_i+S_u)\big(\Gamma(S_k,S_u)\Gamma(S_j,S_k+S_u)\\
	&&\quadd\qquad\quadd\qquad -\Gamma(S_j,S_u)\Gamma(S_k,S_j+S_u)\big)x^{\lam_i+\lam_j+\lam_k+\lam_u},\\	
\end{eqnarray*}
and
\begin{eqnarray*}
	J_2:=[{x^{\lam_j}}\cdot {x^{\lam_k}}]\cdot({x^{\lam_i}}\cdot  x^{\lam_u})
	&=&
		\Gamma(S_i,S_u)[x^{\lam_j},x^{\lam_k}]x^{\lam_i+\lam_u}\\
			&=&\Gamma(S_i,S_u)\big(\Gamma(S_k,S_i+S_u)\Gamma(S_j,S_k+S_u+S_i)\\
			&&-\Gamma(S_j,S_i+S_u)\Gamma(S_k,S_j+S_u+S_i)\big)x^{\lam_i+\lam_j+\lam_k+\lam_u}.\\	
				&=&\Gamma(S_i,S_u)\big(\Gamma(S_k,S_i+S_u)\Gamma(S_j,S_j+S_u)\\
			&&-\Gamma(S_j,S_i+S_u)\Gamma(S_k,S_k+S_u)\big)x^{\lam_i+\lam_j+\lam_k+\lam_u}.\\		
\end{eqnarray*}
 If $u=0$, then  we have
$$\begin{array}{c}
	\Gamma(S_j,S_k+S_u)=\Gamma(S_j,S_k)=0,\quad \Gamma(S_k,S_j+S_u)=\Gamma(S_k,S_j)=0,\\
	\Gamma(S_k,S_i+S_u)=\Gamma(S_k,S_i)=0,\quad
	\Gamma(S_j, S_i+S_u)=\Gamma(S_j,S_i)=0
\end{array}
$$
and so $J_1=J_2=0$.
If $u=i$, then $\Gamma(S_k,S_u)=\Gamma(S_j, S_u)=0$ so $J_1=0$. Also $\Gamma(S_j,S_j+S_u)=\Gamma(S_j,S_k)=0$ and $\Gamma(S_k,S_k+S_u)=\Gamma(S_k,S_{j})=0$, so $J_2=0$. If $i=j$, then
$\Gamma(S_i,S_i+S_u)=\Gamma(S_i,S_k)=0$ and $\Gamma(S_i,S_u)=0$ so $J_1=0=J_2$. The case $u=k$ holds by symmetry.
Finally if $u\not\in\{0,i,j,k\}$, then $\Gamma(S_k, S_u)=\Gamma(S_j,S_u)=\Gamma(S_i,S_u)=0$ and so again $J_1=0=J_2$.\qed

\begin{lem}\label{coia1}
	If $\lam_i,\mu_i\in S_i$ and $\lam_j\in S_j$, then
	$$
	[L_{x^{\mu_i}},[L_{x^{\lam_i}},L_{x^{\lam_j}}]]=\left\{
	\begin{array}{ll}
		0&\hbox{if $i=j$ or $0\in\{i,j\}$},\\		
		L_{x^{\mu_i+\lam_i+\lam_j}}&\hbox{otherwise},
	\end{array}
	\right.
	$$
	and so $[L_{x^{\mu_i}},[L_{x^{\lam_i}},L_{x^{\lam_j}}]]=\ep L_{x^{\mu_i+\lam_i+\lam_j}}$, $\ep\in\{0,1\}$.
\end{lem}

\proof From Lemma \ref{final1},  we have 	$[L_{x^{\mu_i}},[L_{x^{\lam_i}},L_{x^{\lam_j}}]]=-L_{[L_{x^{\lam_i}},L_{x^{\lam_j}}](x^{\mu_i})}$, and this expression is $0$
if either $i=j$ or $0\in\{i,j\}$.
Now assume that $1\leq i\not=j$. Then
\begin{eqnarray*}[L_{x^{\lam_i}},L_{x^{\lam_j}}](x^{\mu_i})
	&=&
	x^{\lam_i}\cdot(\stackrel{0}{\overbrace{x^{\lam_j}\cdot x^{\mu_i}}})-x^{\lam_j}\cdot (x^{\lam_i}\cdot x^{\mu_i})\\
	&=&
	\big(\stackrel{0}{\overbrace{\Gamma(S_j,S_i)}}\Gamma(S_i,S_j+S_i)-\stackrel{0}{\overbrace{\Gamma(S_i,S_i)\Gamma(S_j,S_i+S_i)}}\big)x^{\lam_i+\lam_j+\mu_i}\\
	&=&	-{x^{\mu_i+\lam_i+\lam_j}}(x^\gamma).
\end{eqnarray*}\qed

\section{Dimension of isotropic root spaces}\setcounter{equation}{0}\label{isotropic}
From now on, we assume that $\Lam=\bbbz\sg_1\oplus\bbbz\sg_2$ has rank $2$. Then up to similarity, there exist two possible semilattices 
$$\Lam=S_0\cup S_1\cup S_2\cup S_3\andd S=S_0\cup S_1\cup S_2,$$
where
$$S_i=\tau_i+2\Lam\hbox{ with } \tau_0=0,\;\tau_1=\sg_1,\;\tau_2=\sg_2\hbox{ and }\tau_3=\sg_1+\sg_2.$$ 

\subsection{The case $\jj_S$}\setcounter{equation}{0}
Set
$\jj=\jj_S=\sum_{\lam\in S}\bbbk x^\lam$.
As previously noted, we interpret $L_{x^\lam}=0$ if $\lam\not\in S$.
Continuing with the same terminology and notation as in earlier sections, our primary objection in this section is to determine $\dim E_\sg$ for $0\not=\sg\in R^0$.
The most challenging part of this, according to (\ref{dim}), involves computing 
\begin{equation}\label{tag}
	(\star):=[L_{x^{\sg+\tau}},L_{x^{-\tau}}](x^\gamma)=\stackrel{I_1}{\overbrace{x^{\sg+\tau}\cdot(x^{-\tau}\cdot x^\gamma)}}-\stackrel{I_2}{\overbrace{x^{-\tau}\cdot(x^{\sg+\tau}\cdot x^\gamma)}},
\end{equation}
for $\tau,\gamma\in\Lam$. This will be addressed in the rest of this subsection; see also \cite[Proposition 1.4.2]{ABFP09}.

We recall from Remark \ref{rok5}(ii) that $\dim\Hom_\bbbz(\Lam,\bbbk)=2$ and so $\dim\scd(\fg)^\sg\leq 1$ for each $0\not=\sg\in\Lam$.
Therefore, $\dim D^\sg,\dim C^\sg\leq 1$ for $0\not=\sg\in\Lam$.

(ii) 
	\begin{lem}\label{lemsim1} Let {$\jj=\jj_S$} where $S=S_0\cup S_1\cup S_2$ and $0\not=\sg\in R^0$. Then 
	\begin{eqnarray*}
		E_\sg&=&
		\left\{\begin{array}{ll}
			\bbbk[L_{x^{\sg+\sg_2}},L_{x^{-\sg_2}}]+D^\sg+C^\sg&\hbox{if }\sg\in  S_3\\
			\bbbk L_{x^\sg}+D^\sg+C^\sg&\hbox{if }\sg\in S_{0}\cup S_1\cup S_2.
		\end{array}\right.\\
	\end{eqnarray*}
	In particular, $1\leq \dim E_\sg\leq 3$,
		for any nonzero isotropic root $\sg$.
\end{lem}

\proof Note that $R^0=S+S=S\cup S_3$. Assume first that $\sg\in S_3$ and $\tau\in S_j$. If $\sg+\tau\not\in S$, then $(\star)=0$. Otherwise, we have $j=1$ or $j=2$. Applying Lemma \ref{new09}, we obtain
$[L_{x^{\sg+\tau}},L_{x^{-\tau}} ]=[L_{x^{\sg+\sg_j}},L_{x^{-\sg_j}}]$. If $j=1$ and $\sg=\sg_1+\sg_2+2\lam$, then
applying Lemma \ref{new09} again we find
 $$[L_{x^{\sg+\sg_1}},L_{x^{-\sg_1}}]=[L_{x^{2\sg_1+\sg_2+2\lam}}, L_{x^{-\sg_1}}]=
		[L_{x^{-\sg_2}}, L_{x^{\sg+\sg_2}}].$$
		This concludes the proof for the case $\sg\in S_3$.
		
		Next assume that $\sg\in S$, say $\sg\in S_i$ for some $i=0,1,2$.  Assume also that $\tau\in S_j$ for some $j=0,1,2$.
		Then either $\sg+\tau\not\in S$ or one of $\sg+\tau$, $\tau$ belongs to $S_0$. In either case we have
		$[L_{x^{\sg+\tau}},L_{x^{-\tau}} ]=0$ by Lemma \ref{final1}.
\qed

\begin{lem}\label{lemsim5} Let $S=\Lam$, {$\jj=\jj_S$} and $0\not=\sg\in R^0$. Then
	$$E_\sg=\left\{\begin{array}{ll}
		\bbbk L_{x^\sg}\oplus\bbbk [L_{x^{\sg+\sg_2}},L_{x^{-\sg_2}}]\oplus D^\sg\oplus C^\sg&\hbox{if }\sg\in S_1\cup S_3\\
		\bbbk L_{x^\sg}\oplus\bbbk [L_{x^{\sg+\sg_1}},L_{x^{-\sg_1}}]\oplus D^\sg\oplus C^\sg&\hbox{if }\sg\in S_2\\
		\bbbk L_{x^\sg}\oplus D^\sg\oplus C^\sg&\hbox{if }\sg\in S_0.
	\end{array}\right.	
	$$
	In particular,
	$1\leq \dim E_\sg\leq 4$.
\end{lem}

\proof If either of $\sg$, $\tau$ or $\sg+\tau$ is in $S_0$, then we get $(\star)=0$ by Lemma \ref{final1}(ii).
So we may assume that $\sg,\tau\not\in S_0$ and that $\sg$ and $\tau$ are in different $S_i$.
If $\sg\in S_i$ and $\tau\in S_j$, $1\leq i\not=j\leq 3$, then
by Lemma \ref{new09}, 
$[L_{x^{\sg+\tau}},L_{x^{-\tau}}]=[L_{x^{\sg+\tau_j}},L_{x^{-\tau_j}}]$, $\tau_1=\sg_1$, $\tau_2=\sg_2$.
Now if $\sg=\sg_1+\sg_2+2\lam\in S_3$, then using Lemma \ref{new09},
\begin{eqnarray*}
[L_{x^{\sg+\sg_1}},L_{x^{-\sg_1}}]
&=&
[L_{x^{2\sg_1+\sg_2+2\lam}},L_{x^{-\sg_1}}]\\
&=&
[L_{x^{\sg_2}},L_{x^{-\sg_1+2\sg_1+2\lam}}]\\
&=&
[L_{x^{-\sg_2}},L_{x^{\sg_1+2\sg_2+2\lam}}]\\
&=&
[L_{x^{-\sg_2}},L_{x^{\sg+\sg_2}}].\\
\end{eqnarray*}
If $\sg=\sg_1+2\lam\in S_1$, then for $\tau_2=\sg_2$, we have
$
[L_{x^{\sg+\tau_2}},L_{x^{-\tau_2}}]
=
[L_{x^{\sg+\sg_2}},L_{x^{-\sg_2}}],
$
and for $\tau_3=\sg_1+\sg_2$,
\begin{eqnarray*}
[L_{x^{\sg+\tau_3}},L_{x^{-\tau_3}}]&=&
[L_{x^{\sg+\sg_1+\sg_2}},L_{x^{-\sg_1-\sg_2}}]\\
&=&
[L_{x^{2\sg_1+2\lam+\sg_2}},L_{x^{-\sg_1-\sg_2}}]\\
(\hbox{by Lemma \ref{new09}})&=&
 [L_{x^{-\sg_2}},L_{x^{-\sg_1-\sg_2+2\sg_1+2\sg_2+2\lam}}]\\
&=& -[L_{x^{\sg+\sg_2}},L_{x^{-\sg_2}}].
\end{eqnarray*}
Finally, observe that  $[L_{x^{\sg_1+\sg_2}},L_{x^{-\sg_2}}]$ and $[L_{x^{\sg_2+\sg_1}}, L_{x^{-\sg_1}}]$ are linearly independent, since
the action of the first bracket on $x^{\sg_1}$ is zero, while the action of the second one is none-zero.\qed

\subsection{The case $\jj=\bbbk_{\bq}^+$}\setcounter{equation}{0}

	We begin with a quick review of Jordan tori $\bbbk_\bq^+$ and $H(\bbbk_\be,-)$.
	
	Let $\bq=(q_{ij})\in M_ n(\bbbk)$ be a $( \nu\times \nu)$-matrix satisfying $q_{ii}=1=q_{ij}q_{ji}$ for $1\leq i,j\leq \nu$. The {\it quantum torus} $\bbbk_\bq$ associated with $\bq$ is the unital associative algebra defined by generators $x_1^{\pm1},\ldots,x_ \nu^{\pm1}$ and relations $x_ix_i^{-1}=1=x_i^{-1}x_i$, $x_ix_j=q_{ij}x_jx_i$ for $1\leq i,j\leq \nu$. We have $\bbbk_\bq=\sum_{\lam\in\Lam}\bbbk x^\lam$, where $\Lam=\bbbz^\nu$, and for $\lam=(\lam_1,\ldots,\lam_\nu)$,
	$x^\lam:=x_1^{\lam_1}\cdots x_ \nu^{\lam_\nu}$. Define $\eta:\Lam\times\Lam\rightarrow\bbbk^\times$ by
	$x^\lam\cdot x^{\mu}=\eta(\lam,\mu)x^{\lam+\mu}$. In fact $\eta(\lam,\mu)=\prod_{i<j}q_{ji}^{\mu_j\lam_i}$, which makes $\eta$ a bi-additive map.
Then
	$$[x^\lam,x^{\mu}]=(\eta(\lam,\mu)-\eta(\mu,\lam))x^{\lam+\mu}
	=\eta(\mu,\lam)(f(\lam,\mu)-1)x^{\lam+\mu},$$
	where $f(\lam,\mu)=\eta(\lam,\mu)\eta(\mu,\lam)^{-1}.$ The subgroup
$\rad(f):=\{\lam\mid f(\lam,\mu)=1\hbox{ for all }\mu\in\Lam\}$ of $\Lam$ is called the {\it radical} of $f$.
For $\nu=2$, the quantum torus $\bbbk_\bq$ is uniquely determined by a single parameter $q=q_{21}\in\bbbk^\times$.
The {\it plus algebra} $\bbbk_\bq^+$ of $\bbbk_\bq$ is defined by the same underlying vector space as $\bbbk_\bq$ but with the Jordan product given by $x\cdot y:=\frac{1}{2}(xy+yx).$

	\begin{lem}\label{lemsim3} Let $\jj=\bbbk_\bq^+$, and $0\not=\sg\in R^0$. Then 
$$E_\sg=\left\{\begin{array}{ll}
		\bbbk L_{x^\sg}+D^\sg+C^\sg&\hbox{if }\sg\in\rad(f)\\
	\bbbk L_{x^\sg}\oplus\bbbk (L_{x^\sg}-r_\sg)+D^\sg+C^\sg&\hbox{ if }\sg\not\in\rad(f).
\end{array}\right.
$$
 where $r_\sg:\bbbk_\bq\rightarrow\bbbk_\bq$ is given by $r_\sg(x^\lam)=\eta(\sg,\lam)x^{\sg+\lam}$. In particular,
	$1\leq \dim E_\sg\leq 3 $ if $\sg\in\rad(f)$, and
$1\leq\dim E_\sg\leq 4$ if $\sg\not\in\rad(f).$
\end{lem}

\proof  We fix $0\not= \sg\in R^0=\Lam$. Let $\tau,\gamma\in\Lam$.
Then 
\begin{eqnarray*}
	[L_{x^{\sg+\tau}},L_{x^{-\tau}}](x^\gamma)
	&=&x^{\sg+\tau}\cdot (x^{-\tau}\cdot x^\gamma)-x^{-\tau}\cdot (x^{\sg+\tau}\cdot x^\gamma)\\
	&=&\frac{1}{2}x^{\sg+\tau}\cdot (x^{-\tau}x^\gamma+x^\gamma x^{-\tau})\\
	&&-\frac{1}{2}x^{-\tau}\cdot(x^{\sg+\tau}x^\gamma+x^\gamma x^{\sg+\tau})\\
		&=&\frac{1}{4}\big(x^{\sg+\tau}(x^{-\tau}x^\gamma+x^\gamma x^{-\tau})+
		(x^{-\tau}x^\gamma+x^\gamma x^{-\tau})x^{\sg+\tau}\big)\\
	&&-\frac{1}{4}\big(x^{-\tau}(x^{\sg+\tau}x^\gamma+x^\gamma x^{\sg+\tau})
	+(x^{\sg+\tau}x^\gamma+x^\gamma x^{\sg+\tau})x^{-\tau}\big)\\
	&=&\frac{1}{4}\big(x^{\sg+\tau}x^{-\tau}x^\gamma-x^{-\tau}x^{\sg+\tau}x^\gamma+
	x^\gamma x^{-\tau}x^{\sg+\tau}
	-x^\gamma x^{\sg+\tau}x^{-\tau}\big)\\
	&=&
	\frac{1}{4}(x^{\sg+\tau}x^{-\tau}x^\gamma-x^{-\tau}x^{\sg+\tau}x^\gamma)+
	\frac{1}{4}
	x^\gamma (x^{-\tau}x^{\sg+\tau}
	-x^{\sg+\tau}x^{-\tau})\\
	&=&
	\eta_\tau(\frac{1}{2} (x^\sg x^\gamma+x^\gamma  x^\sg)- x^\gamma x^\sg)=
\eta_{\tau}(L_{x^\sg}-r_{{\sg}})(x^\gamma),
	\end{eqnarray*}
where
\begin{eqnarray*}
\eta_\tau:&=&\frac{1}{2}(	\eta(\sg+\tau,-\tau)-\eta(-\tau,\sg+\tau))\\
&=&\frac{1}{2}(	\eta(\sg,-\tau)-\eta(-\tau,\sg))\\
&=&\frac{1}{2}\eta(-\tau,\sg)(f(\sg,-\tau)-1),
\end{eqnarray*}
and the linear map $r_\sg:\bbbk_\bq\rightarrow\bbbk_\bq$ is defined by
$r_\sg(x^\gamma)=x^\gamma x^\sg=\eta(\gamma,\sg)x^{\sg+\gamma}.$ If $\sg\in\rad(f)$, then $\eta_\tau=0$ for all $\tau$, so 
$$E_\sg=\left\{\begin{array}{ll}
\bbbk L_{x^\sg}\oplus\bbbk (L_{x^\sg}-r_\sg)+D^\sg+C^\sg&\hbox{ if }\sg\not\in\rad(f)\\
\bbbk L_{x^\sg}+D^\sg+C^\sg&\hbox{if }\sg\in\rad(f).
\end{array}\right.
$$
\qed

\subsection{The case $H(\bbbk_\be,-)$}\setcounter{equation}{0}
Consider the elementary quantum matrix $\be=(e_{ij})$ with $e_{ij}=1$ or $-1$ for all $i,j$. Equip the quantum torus $\bbbk_\be$ with the involution $\bar{\;}$ satisfying $\overline{x_i}=x_i$ for all $i$. Then $H(\bbbk_\be,-)=\{x\in\bbbk_\be\mid\overline{x}=x\}$ is a Jordan subalgebra of $\bbbk_\be^+$. Since $\bar{\;}$ preserves the grading of $\bbbk_\be^+$, we see that $H(\bbbk_\be, -)$ is a Jordan $\Lam$-torus  with $H(\bbbk_\be,-)=\sum_{\lam\in\Lam}\left(\bbbk x^\lam\cap H(\bbbk_\be,-)\right)$, see\cite[Example 4.3 (2)]{Yos00}.

	\begin{lem}\label{lemsim4} 
	Let $\jj=H(\bbbk_\be,-)$ and $0\not=\sg\in R^0$. Then 
	$$
	E_\sg=\left\{\begin{array}{ll}
		\bbbk L_{x^\sg}\oplus D^\sg\oplus C^\sg&\hbox{if }-=\id,\hbox{ or }\sg\in S_0\cup S_1\cup S_2,\\
		\bbbk [L_{x^{\sg+\sg_1}}, L_{x^{-\sg_1}}]\oplus D^\sg\oplus C^\sg&\hbox{otherwise}.
	\end{array}\right.
	$$		
	In particular $1\leq \dim E_\sg\leq 3$.
\end{lem}

\proof We have either $\be=\left(\begin{array}{cc}
	1&1\\
	1&1
\end{array}\right),
$
or
$\be=\left(\begin{array}{cc}
	1&-1\\
	-1&1
\end{array}\right).
$
In the first case $\bbbk_\be$ is the algebra of Laurent polynomials in variables $x_1, x_2$,
and $\bar{\;}$ is just the identity map. Then $(\star)=0$ (see (\ref{tag}) and so $E_\sg=\bbbk L_{x^\sg}\oplus D^\sg\oplus C^\sg$.

Next, we consider the second case. Then $x_2x_1=-x_1x_2$ and for $\lam=k_1\sg_1+k_2\sg_2\in \Lam$, $x^\lam\in H(\bbbk_\be,-)$ if and only if $k_1k_2\in 2\bbbz$. Thus $H(\bbbk_\be,-)=\sum_{\lam\in S}\bbbk x^\lam$
with $S=S_0\cup S_1\cup S_2$ with $S_i=\sg_i+2\Lam$ where $\sg_0=0$. We now compute $[L_{x^{\sg+\tau}}, L_{x^{-\tau}}](x^\gamma)$.
Fix $0\not= \sg\in R^0=\Lam$. Let $\tau,\gamma\in S$.
From proof of Lemma \ref{lemsim3}, we have
\begin{eqnarray*}
	[L_{x^{\sg+\tau}},L_{x^{-\tau}}](x^\gamma)
	&=&
	\frac{1}{4}(x^{\sg+\tau}x^{-\tau}x^\gamma-x^{-\tau}x^{\sg+\tau}x^\gamma)+
	\frac{1}{4}
	x^\gamma (x^{-\tau}x^{\sg+\tau}
	-x^{\sg+\tau}x^{-\tau}).
\end{eqnarray*}
Thus
\begin{eqnarray*}
	[L_{x^{\sg+\tau}},L_{x^{-\tau}}](x^\gamma)
&=&
\frac{1}{4}[\big(\eta(-\tau,\gamma)\eta(\sg+\tau,-\tau+\gamma)-\eta(\sg+\tau,\gamma)\eta(-\tau,\sg+\tau+\gamma)\big)\\
&+&\big(\eta(\gamma,-\tau)\eta(\gamma-\tau,\sg+\tau)-\eta(\gamma,\sg+\tau)\eta(\gamma+\sg+\tau,-\tau)\big)]x^{\sg+\gamma}\\
&=&
\frac{1}{4}\big(\eta(-\tau,\gamma)\eta(\sg+\tau,-\tau)\eta(\sg+\tau,\gamma)\\
&&\qquad\qquad\qquad-\eta(\sg+\tau,\gamma)\eta(-\tau,\sg+\tau)\eta(-\tau,\gamma)\big)x^{\sg+\gamma}\\
&+&\frac{1}{4}\big(\eta(\gamma,-\tau)\eta(\gamma,\sg+\tau)\eta(-\tau,\sg+\tau)\\
&&\qquad\qquad\qquad-\eta(\gamma,\sg+\tau)\eta(\gamma,-\tau)\eta(\sg+\tau,-\tau)\big)x^{\sg+\gamma}\\
&=&
\eta_\tau\eta_{\gamma,\tau}
\end{eqnarray*}
where
$$
\eta_\tau:=\frac{1}{2}(	\eta(\sg+\tau,-\tau)-\eta(-\tau,\sg+\tau)),
$$
and
$$\eta_{\gamma,\tau}:=\frac{1}{2}\big(\eta(-\tau,\gamma)\eta(\sg+\tau,\gamma)-\eta(\gamma,-\tau)\eta(\gamma,\sg+\tau)\big).
$$
To compute $\eta_\tau$ and $\eta_{\gamma,\tau}$, we note that if $\lam_i,\lam'_i\in S_i$, $i=0,1,2$, then
$$\eta(\lam_i,\lam_0)=\eta(\lam_0,\lam_i)=\eta(\lam_i,\lam_i)=\eta(\lam_2,\lam_1)=1\andd \eta(\lam_1,\lam_2)=-1.$$
 
Now assume first that $\sg\in S$. If $\sg\in S_0$, then as $\tau\in S$, we have $\sg+\tau,-\tau\in S_i$ for some $i=0,1,2$ and so $\eta_\tau=0$. If $\sg\in S_1$ then as $\sg+\tau\in S$, we have $\tau\in S_0\cup S_1$ and so again $\eta_\tau=0$.
By symmetry, $\eta_\tau=0$ if $\sg\in S_2$. Thus $\eta_\tau\eta_{\gamma,\tau}=0$ if $\sg\in S$.

Assume next that $\sg\in\Lam\setminus S$, namely $\sg=\sg_1+\sg_2+2\lam$ for some $\lam\in\Lam$.
Since $\sg+\tau\in S$, we have $\tau\in S_1\cup S_2$. Now $\sg+\tau\in S_2$ if $\tau\in S_1$ and $\sg+\tau\in S_1$ if $\tau\in S_2$. Therefore
$$\eta_\tau=\left\{\begin{array}{rl}
	1&\hbox{if }\tau\in S_1\\
	-1&\hbox{if }\tau\in S_2.
	\end{array}\right.
	$$
	We now compute $\eta_{\gamma,\tau}$. If either of $\tau$ or $\gamma\in S_0$, then $\eta_{\gamma,\tau}=0$.
	If $\tau,\gamma\in S_1$, then $\sg+\tau\in S_2$ and so $\eta_{\gamma,\tau}=1$. In general with similar computations, we get
	$$\eta_{\gamma,\tau}=\left\{\begin{array}{rl}
	0	&\hbox{if }\tau\hbox{ or }\gamma\in S_0\\
		1&\hbox{if }\tau,\gamma\in S_1\\
			-1&\hbox{if }\tau,\gamma\in S_2\\
				-1&\hbox{if }\tau\in S_1,\gamma\in S_2\\
				1&\hbox{if }\tau\in S_2,\gamma\in S_1.\\
			\end{array}\right.
	$$
	Thus $\eta_{\sg_1}\eta_{\gamma,\sg_1}=-\eta_{\sg_2}\eta_{,\gamma,\sg_2},$ and so $[L_{x^{\sg+\tau}}, L_{x^{-\tau}}](x^\gamma)=\ep [L_{x^{\sg+\sg_1}},L_{x^{-\sg_1}}](x^\gamma),$ $\ep\in\{0,\pm1\}$.
	This completes the proof.\qed
	
Putting together Lemmas \ref{lemsim1}, \ref{lemsim5}, \ref{lemsim3} and \ref{lemsim4}, we have the following result.

	\begin{pro}
	Let $E=E(\fg,D,\kappa)$ be an elliptic extended affine Lie algebra of type $A_1$ where $\fg=\hbox{TKK}(\jj)$. If $0\not=\sg\in R^0$, then
	$\dim E_\sg$ is given by
	$$\left\{\begin{array}{ll}
		2+\dim D^\sg+\dim C^{\sg}&\hbox{if }\jj=J_\Lam\hbox{ and }\sg\not\in2\Lam,\hbox{ or }\jj=\bbbk_\bq^+\hbox{ and }\sg\not\in\rad(f),\\
		1+\dim D^\sg+\dim C^{\sg}&\hbox{otherwise}.
	\end{array}\right.
	$$ 
	In particular,  $2\leq \dim E_\sg\leq 4$ if $\jj=J_\Lam$ and $\sg\not\in 2\Lam$,  or $\jj=\bbbk_\bq^+$ and $\sg\not\in\rad(f),$
		and
		$1\leq \dim E_\sg\leq 3\hbox{ otherwise}.
	$
\end{pro}

\section{Chevalley bases}\setcounter{equation}{0}\label{finalsec}
We continue using notation established in the previous sections. Let $E=E(\fg,D,\kappa)$ be an elliptic extended affine Lie algebra of type $A_1$ with $\kappa=0$. We have $R=(S+S)\cup (\pm\a+S)$, where $S$ is a semilattice with $\la S\ra=\Lam=\bbbz\sg_1\oplus\bbbz\sg_2$.

To discuss the Chevalley bases for $E$, we consider the most general case by assuming $D=\scd(\fg)$. The cases $D\sub\scd(\fg)$ can be derived from this general setting. 
Let $0\not=\sg\in R^0$. We have
$D^\sg=\chi^\sg\{\partial_\theta\in\dd\mid \theta(\sg)=0\}$, where $\dd=\{\partial_\theta\mid\theta\in\Hom_\bbbz(\Lam,\bbbk)$\}. Here $\Hom_\bbbz(\Lam,\bbbk)=\bbbk\theta_1\oplus\bbbk\theta_2$, where $\theta_i (\sg_j)=\delta_{ij}$.
Suppose $0\not=\sg\in\Lam$, say  $\sg=k_1\sg_1+k_2\sg_2$, $k_i\in\bbbz$, and assume without loos of generality that
$k_2\not=0$. Define $\theta_\sg:=k_2\theta_1-k_1\theta_2$. Then 
$$\theta_\sg(\sg)=0,\quad 
\theta_\sg(\Lam)\sub\bbbz,
$$ and thus $$D^\sg=\bbbk\chi^\sg\partial_{\theta_\sg}.$$ 
Define $c^\sg\in (D^\lam)^\star$ by $c^\sg(\chi^\lam\partial_{\theta_\lam})=\delta_{-\sg,\lam}$. Then, by recalling that $((gr D)^\star)^\lam=(D^{-\lam})^\star$, we have $C^\sg=\bbbk c^{\sg}$.

\subsection{Tables of Chevalley bases}
We explicitly introduce a Chevalley basis $\bb$ for $E$ as follows. Recall that $\Lam=S_0\cup S_1\cup S_2\cup S_3$, where
$$S_i=\tau_i+2\Lam\hbox{ with }\tau_0=0,\;\tau_1=\sg_1,\;\tau_2=\sg_2\hbox{ and }\tau_3=\sg_1+\sg_2.$$
For each $\beta\in R$, we assign the elements of {$\bb\cap E_\beta$} by:

\begin{table}[htb]
	\caption{The elements of  {$\bb\cap E_\beta$}}
	\label{tab:fq}
	\vspace{-4mm}
	\scriptsize{
		\begin{tabular}
			[c]{|c |   c |} \hline
						\diagbox[width=10em]{$\qquad\qquad\beta$}{$\jj=\jj_S, S\not=\Lam$}& $\bb\cap E_\beta$  \\
					\whline
					$\a+\lam\in R^\times$ &$x_\beta=\sqrt{2}x^\lam$ \\
					\hline
				
					$-\a+\lam$ &$x_\beta=\sqrt{2}\bar{x^\lam}$  \\
					\hline 
					$0\not=\lam\in S_0\cup S_1\cup S_2$ & $	x^1_\beta=L_{x^\lam},$ $x^2_\beta=\chi^\lam\partial_{\theta_\lam},x^3_\beta=c^\lam$ \\
					\hline
					$\lam\in S_3$ & $x_\beta^1=[L_{x^{\lam+\sg_2}},L_{x^{-\sg_2}}], x^2_\beta=\chi^\lam\partial_{\theta_\lam},x^3_\beta=c^\lam$\\
					\hline
					$0$ & $x_0^1=h_\a, x_0^2=\partial_{\theta_1}, x_0^3=\partial_{\theta_2}, x^4_0=c^1,x^5_0=c^2,$ with $c^i(\partial_{\theta_j})=\delta_{ij}$ \\
					\hline
				\end{tabular}
			}
		\end{table}
	
\begin{table}[htb]
	\caption{The elements of  $\bb\cap E_\beta$}
	\label{tab:fq2}
	\vspace{-4mm}
	\scriptsize{
		\begin{tabular}
			[c]{|c |   c |} \hline
			\diagbox[width=10em]{$\qquad\qquad\beta$}{ $\jj=\jj_\Lam$}& $\bb\cap E_\beta$ \\
			\whline
			$\a+\lam\in R^\times$ &$x_\beta=\sqrt{2}x^\lam$ \\
			\hline
			
			$-\a+\lam$ &$x_\beta=\sqrt{2}\bar{x^\lam}$  \\
			\hline 
			$0\not=\lam\in S_1\cup S_3$ & $	x^1_\beta=L_{x^\lam},
			 x_\beta^2=[L_{x^{\lam+\sg_2}}, L_{x^{-\sg_2}}], x^3_\beta=\chi^\lam\partial_{\theta_\lam},x^4_\beta=c^\lam$ \\
			 \hline 
			$0\not=\lam\in S_2$ & $	x^1_\beta=L_{x^\lam},
			 x_\beta^2=[L_{x^{\lam+\sg_1}}, L_{x^{-\sg_1}}], x^3_\beta=\chi^\lam\partial_{\theta_\lam},x^4_\beta=c^\lam$ \\
			 \hline 
			$0\not=\lam\in S_0$ & $	x^1_\beta=L_{x^\lam},
 x^2_\beta=\chi^\lam\partial_{\theta_\lam},x^3_\beta=c^\lam$ \\
			\hline
					$0$ & $x_0^1=h_\a, x_0^2=\partial_{\theta_1}, x_0^3=\partial_{\theta_2}, x^4_0=c^1,x^5_0=c^2,$ with $c^i(\partial_{\theta_j})=\delta_{ij}$ \\
			\hline
		\end{tabular}
	}
\end{table}

\begin{table}[htb]
	\caption{The elements of  $\bb\cap E_\beta$}
	\label{tab:fq3}
	\vspace{-4mm}
	\scriptsize{
		\begin{tabular}
			[c]{|c |   c |} \hline
			\diagbox[width=10em]{$\qquad\qquad\beta$}{$\jj=\bbbk_\bq^+$}& $ \bb\cap E_\beta$ \\
			\whline
			$\a+\lam\in R^\times$ &$x_\beta=\sqrt{2}x^\lam$ \\
			\hline
			
			$-\a+\lam$ &$x_\beta=\sqrt{2}\bar{x^\lam}$  \\
			\hline 
			$0\not=\lam\not\in\rad(f)$ & $x^1_\beta=L_{x^\lam}, x_\beta^2=L_{x^{\lam}}-r_\lam, x_\beta^3=\chi^\lam\partial_{\theta_\lam},x^4_\beta=c^\lam$ \\
			\hline 
			$0\not=\lam\in\rad(f)$ & $x^1_\beta=L_{x^\lam}, x_\beta^2=\chi^\lam\partial_{\theta_\lam},x^3_\beta=c^\lam$\\
			\hline
			$0$ & $x_0^1=h_\a, x_0^2=\partial_{\theta_1}, x_0^3=\partial_{\theta_2}, x^4_0=c^1,x^5_0=c^2,$ with $c^i(\partial_{\theta_j})=\delta_{ij}$ \\
			\hline
		\end{tabular}
	}
\end{table}

\begin{table}[htb]
	\caption{The elements of  $\bb\cap E_\beta$}
	\label{tab:fq4}
	\vspace{-4mm}
	\scriptsize{
		\begin{tabular}
			[c]{|c |   c |} \hline
			\diagbox[width=10em]{$\qquad\qquad\beta$}{$\jj=H(\bbbk_\be,-)$}& $\bb\cap E_\beta$ \\
			\whline
			$\a+\lam\in R^\times$ &$x_\beta=\sqrt{2}x^\lam$ \\
			\hline
			
			$-\a+\lam$ &$x_\beta=\sqrt{2}\bar{x^\lam}$  \\
			\hline 
			$0\not=\lam\not\in S_3$ & $x^1_\beta=L_{x^\lam}, x_\beta^2=\chi^\lam\partial_{\theta_\lam},x^3_\beta=c^\lam$ \\
			\hline 
			$\lam\in S_3$ & $x^1_\beta=r_{\lam}, x_\beta^2=\chi^\lam\partial_{\theta_\lam},x^3_\beta=c^\lam$\\
			\hline
			$0$ & $x_0^1=h_\a, x_0^2=\partial_{\theta_1}, x_0^3=\partial_{\theta_2}, x^4_0=c^1,x^5_0=c^2,$ with $c^i(\partial_{\theta_j})=\delta_{ij}$. \\
			\hline
		\end{tabular}
	}
\end{table}

\subsection{The main result}
We now proceed to show that the set $\bb$ is indeed a Chevalley basis for the extended affine Lie algebra $E=E(\fg,D,\kappa)$ with $D=\scd(\fg)$ and $\kappa=0$.  Specifically, $\bb$ serves as a $\bbbk$-basis of $E$ and satisfies
$[\bb,\bb]\sub\bbbz\bb.$
We assume the basis $\{\chi^\gamma\mid\gamma\in\Gamma\}$ of $\cc(\fg)$ is normalized so that $\chi^\gamma(x^\lam)=x^{\gamma+\lam}$; see Subsection \ref{const}.

\begin{pro}\label{pro123} Let $\fg=\hbox{TKK}(\jj)$ where $\jj$ is one of rank $2$ Jordan tori $\jj=\jj_S$, $\bbbk^+_\be$ or $H(\bbbk_\be,-).$ Let $D=\scd(\fg)$ and $\kappa=0$. Then the set $\bb$ described in Tables \ref{tab:fq}--\ref{tab:fq4} is a Chevalley basis for the elliptic extended affine Lie algebra $E=E(\fg,D,\kappa)$.
	\end{pro}

\proof Let $\fg$, $D$, $\kappa$ and $\bb$ be as in the statement. By Lemmas \ref{lemsim1},
\ref{lemsim5}, \ref{lemsim3} and \ref{lemsim4}  the set $\bb$ forms a $\bbbk$-basis for $E$. To prove that $\bb$ is a Chevalley basis, we need to show that $[x,y]\in\bbbz\bb$ for all $x,y\in\bb$. We
 proceed with the following computations.
 
 From Lemma \ref{new09}, we see that
 $$
 \begin{array}{c}
 	\hbox{(I)}:=[L_{x^{\lam_j}},L_{x^{\lam_i}}]=[L_{x^{\lam_j+\lam_i-\sg_i}}, L_{x^{-\sg_i}}]\in\bbbz\bb,\vspace{2mm}\\
 	(1\leq i < j\leq 3,\;\lam_i\in S_i,\lam_j\in S_j).
 \end{array}
 $$
 Also
\begin{eqnarray*}
[x_{\a+\lam},x_{-\a+\gamma}]_E &=&[\sqrt{2}x^\lam,\sqrt{2}\overline{x^\gamma}]_E\\
	&=&[\sqrt{2}x^\lam,\sqrt{2}\overline{x^\gamma}]_\fg		+c_D(\sqrt{2} x^\lam,\sqrt{2}\overline{x^\gamma})\\	
		 &=&2L_{x^{\lam}\cdot {x^\gamma}}+2[L_{x^{\lam}},L_{x^\gamma}]+2c_D(x^\lam,\overline{x^\gamma})\\
		 &\in&\ep L_{x^{\lam+\gamma}}+2[L_{x^{\lam}},L_{x^\gamma}]+2c_D(x^\lam,\overline{x^\gamma}),
 \end{eqnarray*}
 $\ep\in\{0,\pm1\}$.
Now for $d=\chi^\mu\partial_{\theta_\mu}\in D^\mu$, $\mu\not=0$, we have
\begin{eqnarray*}
	c_D(x^\lam,\overline{x^\gamma})(\chi^\mu\partial_{\theta_\mu})&=&
	(\chi^\mu\partial_{\theta_\mu}(x^\lam),\overline{x^\gamma})\\
	&=&
	\theta_\mu(\lam)\left\{\begin{array}{ll}
		1&\hbox{if }\mu+\lam=-\gamma\\
		0&\hbox{otherwise}.
		\end{array}\right.\\
	&=&
		\theta_{-\gamma}(\lam)\left\{\begin{array}{ll}
		1&\hbox{if }\mu+\lam=-\gamma\\
		0&\hbox{otherwise}.
	\end{array}\right.\\
		&=&
			\theta_{-\gamma}(\lam)c^{\lam+\gamma}(\chi^\mu\partial_{\theta_\mu}),	
	\end{eqnarray*}
	and
	\begin{eqnarray*}
			c_D(x^\lam,\overline{x^\gamma})(\partial_{\theta_i})&=&
		(\partial_{\theta_i}(x^\lam),\overline{x^\gamma})\\
		&=&
		\theta_i(\lam)\left\{\begin{array}{ll}
			1&\hbox{if }\lam=-\gamma\\
			0&\hbox{otherwise}.
		\end{array}
		\right.\\
		&=&
		\delta_{\lam,-\gamma}(\theta_1(\lam)c^{1}+\theta_2(\lam)c^2)(\partial_{\theta_i}).		
\end{eqnarray*}		
which imply $c_D(x^\lam,\overline{x^\gamma})\in\bbbz\bb$. Thus by (I),

$$\hbox{(II)}:=[x_{\a+\lam},x_{-\a+\gamma}]\in\bbbz\bb.$$
Also,
\begin{eqnarray*}\hbox{(III)}:=[\chi^\mu\partial_{\theta_\mu},\sqrt{2}x^\lam]
	=
	\sqrt{2}\chi^\mu\partial_{\theta_\mu}(x^\lam)
	=
	\theta_\mu(\lam)\sqrt{2}x^{\mu+\lam} \in\bbbz\bb,
	\end{eqnarray*}
$$\hbox{(IV)}:=[\partial_{\theta_i},\sqrt{2}x^\lam]=\theta_i(\lam)\sqrt{2}x^\lam\in\bbbz\bb,$$
\begin{eqnarray*}
\hbox{(V)}:=[\chi^\mu\partial_{\theta_\mu},\chi^\nu\partial_{\theta_\nu}]&=&
\chi^{\mu+\nu}(\theta_\mu(\nu)\partial_{\theta_\nu}-\theta_\nu(\mu)\partial_{\theta_\mu})\\
&=&
\chi^{\mu+\nu}(\theta_\mu(\nu)\partial_{\theta_\nu}+\theta_\mu(\nu)\partial_{\theta_\mu})\\
&=&
\theta_\mu(\nu)\chi^{\mu+\nu}\theta_{\mu+\nu}\in\bbbz\bb.
\end{eqnarray*}

Next, we note that
$
[L_{x^\sg},\sqrt{2}x^\lam]=
L_{x^\sg}(\sqrt{2}x^{\lam})=\sqrt{2}x^\sg\cdot x^\lam.
$
Now if $\jj\not=\bbbk^+_\bq$, then $x^\sg\cdot x^\lam\in\ep x^{\sg+\lam}$, $\ep\in\{0,\pm1\}$. If $\jj=\bbbk_\bq^+$ with $\bq$ not elementary, then we can find $\sg,\lam$ such that $\eta(\sg,\lam)\not\in\bbbz$. If $\bq$ is elementary, clearly
$\eta(\sg,\lam)\in\bbbz$ for all $\sg,\lam$.
Thus
$$\begin{array}{c}
	\hbox{(VI)}:=[L_{x^\sg},\sqrt{2}x^\lam]\in\bbbz\sqrt{2}x^{\sg+\lam}\in\bbbz\bb,\vspace{2mm}
	\\
(\jj\not=\bbbk^+_\bq,\,\bq=\hbox{non-elementary}.)\\ 
\end{array}
$$

Next, we have
\begin{eqnarray*}
	\hbox{(VII)}:=[[L_{x^{\lam+\sg_2}},L_{x^{-\sg_2}}],L_{x^\mu}]=
	L_{[L_{x^{\lam+\sg_2}},L_{x^{-\sg_2}}](x^\mu)}
	\in\bbbz L_{x^{\lam+\mu}}\in\bbbz\bb,
\end{eqnarray*}
so using this  and the Jacobi identity, we get
$$
\begin{array}{c}
	\hbox{(VIII)}:=[[L_{x^{\lam+\sg_i}},L_{x^{-\sg_i}}],[L_{x^{\lam'+\sg_j}},L_{x^{-\sg_j}}]]\in\bbbz\bb,\vspace{2mm}
\\
(i,j=1,2).
	\end{array}$$

Next, since $\chi^\mu\partial_\theta$ is a derivation, we have
$$\hbox{(IX)}:=[\chi^\mu\partial_{\theta_\mu},L_{x^\lam}]=L_{\chi^\mu\partial_{\theta_\mu}(x^\lam)}=\theta_\mu(\lam)L_{\chi^\mu(x^\lam)}\in\bbbz\bb,$$
and
\begin{eqnarray*}
	 \hbox{(X)}:=[\chi^\mu\partial_{\theta_\mu},[L_{x^{\lam+\sg_2}},L_{x^{-\sg_2}}]]
	 &=&-[L_{x^{-\sg_2}},[\chi^\mu\partial_{\theta_\mu},L_{x^{\lam+\sg_2}}]]\\
&&	 \qquad\qquad+	[L_{x^{\lam+\sg_2}},[\chi^\mu\partial_{\theta_\mu},L_{x^{-\sg_2}}]]\\
	 &=&
	  -\theta_\mu(\lam+\sg_2)[L_{x^{-\sg_2}},L_{\chi^\mu (x^{\lam+\sg_2})}]\\
	&& \qquad\qquad-	\theta_\mu(\sg_2)[L_{x^{\lam+\sg_2}},L_{\chi^{\mu}(x^{-\sg_2})}]\\
(\hbox{using (I)})&\in&\bbbz\bb.
	 \end{eqnarray*}

Next, we note that
\begin{eqnarray*}	
	[r_\lam, L_{x^\mu}](x^\gamma)
	&=&
	r_\lam L_{x^\mu}(x^\gamma)-L_{x^\mu} r_\lam(x^\gamma)\\
	&=&
	\frac{1}{2}r_\lam (x^\mu x^\gamma+x^\gamma x^\mu)-\frac{1}{2}\eta(\lam,\gamma) (x^\mu x^{\lam+\gamma}+x^{\lam+\gamma} x^\mu)\\
&=&
	\frac{1}{2}(\eta(\mu,\gamma)+\eta(\gamma,\mu))r_\lam(x^{\gamma+\mu})\\
	&&\qquad
	-\frac{1}{2}\eta(\lam,\gamma)(\eta(\mu,\lam+\gamma)+\eta(\lam+\gamma,\mu))x^{\lam+\mu+\gamma}\\
	&=&
		\frac{1}{2}(\eta(\mu,\gamma)+\eta(\gamma,\mu))\eta(\lam,\gamma+\mu)x^{\gamma+\mu+\lam}\\
		&&\qquad
		-\frac{1}{2}\eta(\lam,\gamma)(\eta(\mu,\lam+\gamma)+\eta(\lam+\gamma,\mu))x^{\lam+\mu+\gamma}\\	
			&=&\frac{1}{2}[\eta(\mu,\gamma)\eta(\lam,\mu)\eta(\lam,\gamma)+\eta(\gamma,\mu)\eta(\lam,\gamma)\eta(\lam,\mu)\\
&&\qquad		- (\eta(\lam,\gamma)\eta(\mu,\lam)\eta(\mu,\gamma)
		-\eta(\lam,\gamma)\eta(\lam,\mu)\eta(\gamma,\mu)]x^{\lam+\mu+\gamma}\\
			&=&\frac{1}{2}\eta(\mu+\lam,\gamma)(\eta(\lam,\mu)-\eta(\mu,\lam))x^{\lam+\mu+\gamma}\\
			&=&\frac{1}{2}(\eta(\lam,\mu)-\eta(\mu,\lam))x^{\lam+\mu}x^\gamma\\		
				&=&\frac{1}{2}(\eta(\lam,\mu)-\eta(\mu,\lam))x^{\lam+\mu}x^\gamma	+\frac{1}{2}(\eta(\lam,\mu)-\eta(\lam,\mu))x^\gamma x^{\lam+\mu}\\	
				&=&\eta(\lam,\mu)L_{x^{\lam+\mu}}(x^\gamma)
				-(\eta(\mu,\lam)+\eta(\lam,\mu))L_{x^{\lam+\mu}}(x^\gamma)\\
				&=&-	\eta(\mu,\lam)L_{x^{\lam+\mu}}(x^\gamma)\\								
									\end{eqnarray*}	
Thus
$$
\begin{array}{c}\hbox{(XI)}:=	[r_\lam, L_{x^\mu}]=-	\eta(\mu,\lam)L_{x^{\lam+\mu}}\in\bbbz\bb,\vspace{2mm}\\
(\hbox{if $\bq$ is elementary}).
\end{array}
$$	
Also,
\begin{eqnarray*}
	[L_{x^\mu}-r_\mu,\sqrt{2}x^\lam]&=&(L_{x^\mu}-r_\mu)(\sqrt{2}x^\lam)\\
	&=&\frac{\sqrt{2}}{2}(x^\mu\cdot x^\lam-\eta(\mu,\lam)x^{\mu+\lam})
\\
&=&	
\frac{\sqrt{2}}{2}(\eta(\lam,\mu)-\eta(\mu,\lam))x^{\lam+\mu},													\end{eqnarray*}
	and therefore,
	$$
	\begin{array}{c}
		\hbox{(XII)}:=	[L_{x^\mu}-r_\mu,\sqrt{2}x^\lam]\in\bbbz\sqrt{2}x^{\mu+\lam}\in\bbbz\bb,\vspace{2mm}\\
		(\hbox{if }\bq \hbox{ is elementary}).
		\end{array}
	$$	
	
Next, we want to investigate brackets of the form $[\chi^\mu\partial_{\theta_\mu},L_{x^\lam}-r_\lam].$
Since we are interested in elementary case, namely $\bbbk^+_\be$, we have $\eta(\lam,\gamma)\in\{\pm 1\}$ for all
$\lam,\gamma$. Therefore, by the proof of Lemma \ref{lemsim3} we are done if we show that  $[\chi^\mu\partial_{\theta_\mu},[L_{x^{\lam+\tau}},L_{x^{-\tau}}]]\in\bbbz\bb$ for each $\tau$. Now, a similar computation as in (X), gives
\begin{eqnarray*}
	[\chi^\mu\partial_{\theta_\mu},[L_{x^{\lam+\tau}},L_{x^{-\tau}}]]
	&=&
	-\theta_\mu(\lam+\tau)[L_{x^{-\tau}},L_{x^{\lam+\tau+\mu}}]\\
	&& \qquad\qquad-	\theta_\mu(\tau)[L_{x^{\lam+\tau}},L_{x^{-\tau+\mu}}]\\
(\ep\in\{0,\pm1\}, \tau'=\tau-\mu)		&=&
	\ep\theta_\mu(\lam+\tau)(L_{x^{\lam+\mu}}-r_{\lam+\mu})\\
	&& \qquad\qquad-	\theta_\mu(\tau)[L_{x^{\lam+\mu+\tau'}},L_{x^{-\tau'}}]\\
(\ep'\in\{0,\pm 1\})		&=&
	\ep\theta_\mu(\lam+\tau)(L_{x^{\lam+\mu}}-r_{\lam+\mu})\\
		&& \qquad\qquad+\ep'	\theta_\mu(\tau)(L_{x^{\lam+\mu}}-r_{\lam+\mu}).\\
\end{eqnarray*}
Therefore, since $\ep,\ep',\theta_\mu(\lam+\tau),\theta_\mu(\tau)\in\bbbz$, we get
$$
\hbox{(XIII)}:=[\chi^\mu\partial_{\theta_\mu},L_{x^{\lam}}-r_\lam]\in\bbbz\bb.
	$$

Putting together relations (I)-(XIII), we conclude that $[\bb,\bb]\sub\bbbz\bb$, as required.\qed


\begin{bibdiv}
	\begin{biblist}
		

		\bib{AABGP97}{article}{
			label={AABGP97}
			author={{Allison}, Bruce},
			author={{Azam}, Saeid},
			author={{Berman}, Stephen},
			author={{Gao}, Yun},
			author={{Pianzola}, Arturo},
			title={{Extended affine Lie algebras and their root systems}},
			date={1997},
			ISSN={0065-9266; 1947-6221/e},
			journal={{Mem. Am. Math. Soc.}},
			volume={603},
			pages={122},
		}
	


\bib{ABFP09}{article}{
		author={{Allison}, Bruce},
		author={{Berman}, Stephen},
		author={{Faulkner}, J.},
			author={{Pianzola}, Arturo}, 
				Title={Multiloop realization of extended affine Lie algebras
and Lie tori},
journal={Trans. Am. Math. Soc.}
volume={361},
Number={9},
pages={4807-–4842},
year={2009}
}
		
		
		\bib{AG01}{article}{
			author={{Allison}, Bruce},
			author={{Gao}, Yun},
			title={{The root system and the core of an extended affine Lie
					algebra}},
			date={2001},
			ISSN={1022-1824; 1420-9020/e},
			journal={{Sel. Math., New Ser.}},
			volume={7},
			number={2},
			pages={149\ndash 212},
		}
		
		\bib{AFI22}{article}{
	label={AFI22}
	author={Azam, Saeid},
	author={Farahmand~Parsa, Amir},
	author={Izadi~Farhadi, Mehdi},
	title={Integral structures in extended affine {Lie} algebras},
	date={2022},
	ISSN={0021-8693},
	journal={J. Algebra},
	volume={597},
	pages={116\ndash 161},
}
		\bib{AI23}{article}{
	author={{Azam}, Saeid},
	author={{Izadi Farhadi}, Mehdi},
	title={Chevalley involutions for Lie tori and extended affine Lie algebras},
	year={2023},
	journal={J. Algebra},
	Volume={634},
	pages={1--43},
}
	
	
	


	


%


	
	
	
	
	
 \bib{Az06}{article}{
 Author = {Azam, Saeid},
 Title = {Generalized reductive {Lie} algebras: connections with extended affine {Lie} algebras and {Lie} tori},
 Journal = {Can. J. Math.},
 Volume = {58},
 Number = {2},
 Pages = {225--248},
 Year = {2006},
}	

\bib{Az97}{article}{
 author = {Azam, Saeid},
 title = {Nonreduced extended affine root systems of nullity {{\(3\)}}},
 journal = {Commun. Algebra},
 volume = {25},
 number = {11},
 pages = {3617--3654},
 year = {1997},
}

		
		\bib{BGK96}{article}{
			author={Berman, Stephen},
			author={Gao, Yun},
			author={Krylyuk, Yaroslav~S.},
			title={Quantum tori and the structure of elliptic quasi-simple {Lie}
				algebras},
			date={1996},
			ISSN={0022-1236},
			journal={J. Funct. Anal.},
			volume={135},
			number={2},
			pages={339\ndash 389},
			url={semanticscholar.org/paper/ec0c3961f4e077104c5d099700fb975dd7aa01e4},
		}
		
		\bib{BGKN95}{article}{
			author={Berman, Stephen},
			author={Gao, Yun},
			author={Krylyuk, Yaroslav},
			author={Neher, Erhard},
			title={The alternative torus and the structure of elliptic quasi-simple
				{Lie} algebras of type {{\(A_ 2\)}}},
			date={1995},
			ISSN={0002-9947},
			journal={Trans. Am. Math. Soc.},
			volume={347},
			number={11},
			pages={4315\ndash 4363},
		}
%
%
	
\bib{Bou08}{book}{
		Author = {Bourbaki, Nicolas},
		Title = {Elements of mathematics. {Lie} groups and {Lie} algebras. {Chapters} 7--9. {Transl}. from the {French} by {Andrew} {Pressley}},
		Edition = {Paperback reprint of the hardback edition 2005},
		Year = {2008},
		Publisher = {Berlin: Springer},
	}




\bib{Che55}{article}{
Author = {Chevalley, Claude},
Title = {Sur certains groupes simples},
Journal = {T{\^o}hoku Math. J. (2)},
Volume = {7},
Pages = {14--66},
Year = {1955},
}

%
	
		
		\bib{Gar78}{article}{
			author={Garland, Howard},
			title={The arithmetic theory of loop algebras},
			date={1978},
			ISSN={0021-8693},
			journal={J. Algebra},
			volume={53},
			pages={480\ndash 551},
		}
		
		\bib{Gar80}{article}{
			author={Garland, Howard},
			title={The arithmetic theory of loop groups},
			date={1980},
			ISSN={0073-8301},
			journal={Publ. Math., Inst. Hautes {\'E}tud. Sci.},
			volume={52},
			pages={5\ndash 136},
		}
	
	\bib{H-KT90}{article}{
		Author = {H{\o}egh-Krohn, Raphael}
		author={Torresani, Bruno},
		Title = {Classification and construction of quasisimple {Lie} algebras},
		Journal = {J. Funct. Anal.},
		Volume = {89},
		Number = {1},
		Pages = {106--136},
		Year = {1990},
	}
		
		\bib{Hum72}{book}{
			author={{Humphreys}, J.~E.},
			title={{Introduction to Lie algebras and representation theory}},
			publisher={Springer, New York, NY},
			date={1972},
			volume={9},
		}
%
%
%
%
%
%
%
%
	
%
		
		\bib{Mit85}{book}{
			author={{Mitzman}, David},
			title={{Integral bases for affine Lie algebras and their universal
					enveloping algebras}},
			publisher={Contemporary Mathematics, Vol. 40, American Mathematical Society
				(AMS), Providence, RI},
			date={1985},
		}
		
		
		\bib{Neh04}{article}{
			author={{Neher}, Erhard},
			title={{Extended affine Lie algebras}},
			date={2004},
			ISSN={0706-1994},
			journal={{C. R. Math. Acad. Sci., Soc. R. Can.}},
			volume={26},
			number={3},
			pages={90\ndash 96},
		}
		
		
		\bib{Neh11}{incollection}{
			author={{Neher}, Erhard},
			title={{Extended affine Lie algebras and other generalizations of affine
					Lie algebras -- a survey}},
			date={2011},
			booktitle={{Developments and trends in infinite-dimensional Lie theory}},
			publisher={Basel: Birkh\"auser},
			pages={53\ndash 126},
		}
	
		
%
		
		\bib{Ste16}{book}{
				label={Ste16}
		author = {Steinberg, Robert},
		title = {Lectures on {Chevalley} groups},
		series = {Univ. Lect. Ser.},
		volume = {66},
		year = {2016},
		publisher = {Providence, RI: American Mathematical Society (AMS)},
	}
%
%
%
		
		
		\bib{Yos00}{article}{
			author={Yoshii, Yoji},
			title={Coordinate algebras of extended affine {Lie} algebras of type
				{{\(A_1\)}}},
			date={2000},
			ISSN={0021-8693},
			journal={J. Algebra},
			volume={234},
			number={1},
			pages={128\ndash 168},
		}
%
		
%
		
	\end{biblist}
\end{bibdiv}
\end{document}